\documentclass[11pt]{article}

\usepackage{latexsym}
\usepackage{amsfonts}
\usepackage{amsmath}

\usepackage[dvips]{color}

\newtheorem{thrm}{Theorem}[section]
\newtheorem{lemma}[thrm]{Lemma}
\newtheorem{prop}[thrm]{Proposition}

\newtheorem{remark}[thrm]{Remark}

\numberwithin{equation}{section}

\usepackage{enumerate}

\usepackage{amssymb}
\usepackage{mathtools}
\mathtoolsset{showonlyrefs}
\usepackage{mathrsfs}
\usepackage{comment}
\usepackage{hyperref}
\usepackage{xcolor}

\def\E{\mathbb{E} }
\def\P{\mathbb{P} }

\def\R{\mathbb{R} }
\def\N{\mathbb{N} }

\makeatletter
\begin{document}
\allowdisplaybreaks

\title{\Large \bf{Asymptotic expansion for
branching killed Brownian motion with drift}\footnote{The research of this project is supported
     by the National Key R\&D Program of China (No. 2020YFA0712900).}}
\author{ \bf  Haojie Hou \hspace{1mm}\hspace{1mm}
Yan-Xia Ren\footnote{The research of this author is supported by NSFC
(Grant Nos.  12071011 and 12231002) and The Fundamental Research Funds for the Central Universities, Peking University LMEQF.
 } \hspace{1mm}\hspace{1mm} and \hspace{1mm}\hspace{1mm}
Renming Song\thanks{Research supported in part by a grant from the Simons
Foundation
(\#960480, Renming Song).}
\hspace{1mm} }
\date{}
\maketitle

\begin{abstract}

Let $Z_t^{(0,\infty)}$ be the point process  formed by the positions of all particles alive at time $t$ in a branching Brownian motion with drift and killed upon reaching 0.
We study the asymptotic expansions of $Z_t^{(0,\infty)}(A)$ for $A= (a,b)$ and $A=(a,\infty)$ under the assumption that $\sum_{k=1}^\infty k(\log k)^{1+\lambda} p_k <\infty$ for large $\lambda$ in the regime of $\theta \in [0,\sqrt{2})$. These results extend and sharpen the results of
Louidor and Saglietti [J. Stat. Phys, 2020]  and Kesten [Stochastic Process. Appl., 1978].

\end{abstract}

\medskip


\medskip


 \section{Introduction and main results}

\subsection{Introduction}

A branching Brownian motion with drift is a continuous-time Markov process defined as follows: At time $0$, there is a particle at site $x$ and it moves according to a Brownian motion with drift $-\theta$, where $\theta\in \R$. After an exponential time of parameter $1$ independent of the movement,
this particle dies and it splits into $k$ offspring with probability $p_k$.
Each of the offspring independently repeats its parent's behavior from their birth-place. This procedure goes on.
We use $\P_x$ and $\E_x$ to denote the law of this process and the corresponding expectation operator. Let $N(t)$ be the set of particles alive at time $t$ and for each $u\in N(t)$, we will use $X_u(t)$ to denote the position of the particle. For $s<t$ and $u\in N(t)$, we will also use $X_u(s)$ to denote the position of the ancestor of $u$ at time $s$.

Suppose now that $x>0$ and that, once a particle hits $(-\infty, 0]$, we remove it (along with all its possible descendants) from the system.
The resulting branching system is called a branching killed Brownian motion with drift. Let $Z_t^{(0,\infty)}$ denote the point process  formed by the positions of all the particles alive at time $t$ in the
branching killed Brownian motion with drift,
i.e.,
\[
Z_t^{(0,\infty)}:= \sum_{u\in N(t)} 1_{\{\min_{s\leq t} X_u(s)> 0\}}\delta_{X_u(t)}.
\]
Assume that
\[
\sum_{k=0}^\infty kp_k=2.
\]
Kesten \cite{Kesten78} proved that for any $\theta\in \R$,  there exists a constant $C=C(x, \theta)>0$ such that
 \begin{equation}\label{Principal-Engenfunction-2}
 	 \E_x\left(Z_t^{(0,\infty)}((0,\infty))\right) \sim  \left\{\begin{array}{ll}\displaystyle
	 Ct^{-3/2} e^{(1-\frac{\theta^2}{2})t},\quad &\theta>0,\\
	 \displaystyle C t^{-1/2} e^t,\quad&   \theta=0,\\
 	 \displaystyle  C e^t,\quad&\theta<0.
 	\end{array}\right.
 \end{equation}
Consequently, when $\theta \geq \sqrt{2}$, the branching killed Brownian motion with drift will die out.
It was proved in \cite{Kesten78} that, when $\theta <\sqrt{2}$ and $\sum_{k=1}^\infty k^2 p_k< \infty$, the branching killed Brownian motion with drift will survive with positive probability.
In \cite{Kesten78}, Kesten also stated, without proof,
that under the assumption $\sum_{k=0}^\infty k^2 p_k<\infty$, there exists a random variable $W(\theta)$ such that
\begin{align}\label{Extinction-probability}
\P_x \left( W(\theta)>0 \big| Z_t^{(0,\infty)}((0,\infty))>0,\ \forall t>0 \right)=1
\end{align}
and that

(i) If $\theta \in [0, \sqrt{2})$, then $\P_x$-a.s.,  simultaneously for all
intervals $\Delta \subset (0,\infty)$ (finite or infinite), it holds that
\begin{align}\label{SLLN-Non-negative-parameter}
	 \frac{Z_t^{(0,\infty)}(\Delta)}{\E_x \left(Z_t^{(0,\infty)}(\Delta) \right)} \stackrel{t\to\infty}{\longrightarrow} W(\theta).
\end{align}

(ii) If $\theta<0$,  then $\P_x$-a.s.,
\begin{align}\label{SLLN-negative-parameter}
	e^{-t} Z_t^{(0,\infty)}((0,\infty)) \stackrel{t\to\infty}{\longrightarrow} W(\theta).
\end{align}
In \cite{LS20}, Louidor and Saglietti  proved that \eqref{Extinction-probability} and \eqref{SLLN-Non-negative-parameter} hold for the case $\theta \in(0,\sqrt{2})$.

The purpose of this paper is to extend and sharp the main result of \cite{LS20}: We weaken the moment condition from $\sum_{k=1}^\infty k^2p_k<\infty$ to $\sum_{k=1}^\infty k\log^{1+\lambda} k p_k <\infty$ for some $\lambda>0$ and, for any $a\ge 0$, give asymptotic expansions of arbitrary order
for
\[
\frac{Z_{t}^{(0,\infty)} ((a,\infty))}{t^{-3/2}e^{(1-\frac{\theta^2}{2})t}}
\]
under this weaker assumption.
We emphasize here that we do not use the results of  \cite{LS20} in this paper and that, as a consequence of Theorem \ref{thm1} below, we give another proof of \eqref{SLLN-Non-negative-parameter} under the weaker condition.
 We also include the case $\theta=0$.
It is natural to study similar problems for the case $\theta<0$.
We believe things might be different in this case and we plan to tackle this in a future work.

For other recent results on branching killed Brownian motion with drift, see
\cite{BBS15, BBS14, BBS11, BBS13, BBHM17, HH07, HHK06, Liu21, LS23, MS20}
and references therein.
For asymptotic expansions for branching random walks and branching Wiener processes, see
 \cite{Chen01, GL16-1, GL16-2, GL18, GL21, GK17, RRS05}
and references therein.

\subsection{Main results}

We will assume that
\begin{align}\label{LlogL}
	\sum_{k=1}^\infty k(\log k)^{1+\lambda}p_k < \infty,
\end{align}
for some $\lambda>0$.
Let $H_k$ be the $k$-th order Hermite polynomial: $H_0(x):=1$ and for $k\geq 1$,
\[
H_k(x):= \sum_{j=0}^{[k/2]} \frac{k! (-1)^j}{2^j j! (k-2j)!} x^{k-2j}.
\]
It is well known that, if $\{ (B_t)_{t\geq 0}, \Pi_0\}$ is a standard Brownian motion, then,
for any $k\geq 1$,
 $\{t^{k/2} H_k(B_t/\sqrt{t}), \sigma(B_s: s\leq t),\ \Pi_0\}$ is a martingale.
Throughout this paper $\N=\{0,1,2\cdots\}$.
Now for $\theta\in [0,\sqrt{2})$ and $k \in \N$, we define
\begin{align}\label{Martingale}
	M_t^{(2k+1,\theta)}:= e^{-(1-\frac{\theta^2}{2})t}\sum_{u\in N(t)}  1_{\{\min_{s\leq t} X_u(s)> 0\}}e^{\theta X_u(t)} t^{(2k+1)/2} H_{2k+1}\left(\frac{X_u(t)}{\sqrt{t}}\right),
\quad t\geq 0.
\end{align}
We will prove later
(see Proposition \ref{Convergence-rate-martingale}) that, for
any $x>0$, $k\in \N$ and
$\theta \in [0,\sqrt{2})$,
$M_t^{(2k+1,\theta)}$ is a martingale, and if
 \eqref{LlogL} holds for some large $\lambda$, then
$M_t^{(2k+1,\theta)}$ converges to a limit $M_\infty^{(2k+1,\theta)}$ $\P_x$-almost surely and in $L^1(\P_x)$.

\begin{thrm}\label{thm1}
 Assume the drift $\theta$ is in $(0,\sqrt{2})$.
 For any given $m\in \N$, if \eqref{LlogL} holds for  some
   $\lambda >6m +6$,
 then for any $x>0, a\geq 0 $, $\P_x$-almost surely, as $t\to\infty$,
  \begin{align}
  	&\frac{Z_{t}^{(0,\infty)} ((a,\infty))}{t^{-3/2}e^{(1-\frac{\theta^2}{2})t}}\nonumber\\
  	&=- \sqrt{\frac{2}{\pi}}\sum_{\ell=0}^m \frac{H_{2\ell+2}(0)}{t^{\ell}} \sum_{k=0}^{\ell} \frac{M_\infty^{(2k+1,\theta)} }{(2k+1)! (2\ell -2k+1)!}
  	\int_a^\infty z^{2\ell -2k+1} e^{-\theta z} \mathrm{d} z+o(t^{-m}).
  \end{align}
\end{thrm}

 Note that, from
Theorem \ref{thm1}, one can immediately get the asymptotic expansion of
\[
\frac{Z_{t}^{(0,\infty)} ((a,b))}{t^{-3/2}e^{t}}
\]
for any finite interval $(a, b)\subset (0, \infty)$.
For the case $\theta=0$, the result is a little bit different.
For finite intervals, the normalization function is the same as in Theorem \ref{thm1}. For infinite intervals, the
normalization function is different.

\begin{thrm}\label{thrm2}
Assume that the drift $\theta$ is 0.
For any given $m\in \N$, if \eqref{LlogL} holds for  some
  $\lambda >6m +6$, then
	
	(i)  for any $x>0, b>a\geq 0 $, $\P_x$-almost surely, as $t\to\infty$,
	\begin{align}\label{step_44}
		&\frac{Z_{t}^{(0,\infty)} ((a,b))}{t^{-3/2}e^{t}}\nonumber\\
		&=- \sqrt{\frac{2}{\pi}}\sum_{\ell=0}^m \frac{H_{2\ell+2}(0)}{t^{\ell}} \sum_{k=0}^{\ell} \frac{M_\infty^{(2k+1,0)} }{(2k+1)! (2\ell -2k+1)!}
		\int_a^b z^{2\ell -2k+1}  \mathrm{d} z+o(t^{-m});
	\end{align}
	
   (ii) for any $x>0, a\geq 0 $, $\P_x$-almost surely, as $t\to\infty$,
 \begin{align}\label{step_45}
 	&\frac{Z_{t}^{(0,\infty)} ((a,\infty))}{t^{-1/2}e^{t}}= \sqrt{\frac{2}{\pi}}\sum_{\ell=0}^m \frac{H_{2\ell}(0)}{t^{\ell}} \sum_{k=0}^{\ell} \frac{M_\infty^{(2k+1,0)} }{(2k+1)! (2\ell -2k)!} a^{2\ell-2k} +o(t^{-m}).
 \end{align}
\end{thrm}

\begin{remark}
 Note that we only dealt with the case that the branching rate is 1 and the mean number of offspring is 2 in the two theorems above. In the general case when the branching rate is $\beta>0$ and  the mean number of offspring is $\mu>1$, one can use the same argument to prove the following counterpart of Theorem \ref{thm1}:
	Let $\theta\in (0,\sqrt{2\beta(\mu-1)})$. For any given $m\in \N$, if \eqref{LlogL} holds for  some
	  $\lambda >6m +6$,
	then for any $x>0, a\geq 0 $, $\P_x$-almost surely, as $t\to\infty$,
	\begin{align}
		&\frac{Z_{t}^{(0,\infty)} ((a,\infty))}{t^{-3/2}e^{(\beta(\mu-1)-\frac{\theta^2}{2})t}}\nonumber\\
		&=- \sqrt{\frac{2}{\pi}}\sum_{\ell=0}^m \frac{H_{2\ell+2}(0)}{t^{\ell}} \sum_{k=0}^{\ell} \frac{M_\infty^{(2k+1,\theta)} }{(2k+1)! (2\ell -2k+1)!}
		\int_a^\infty z^{2\ell -2k+1} e^{-\theta z} \mathrm{d} z+o(t^{-m}),
	\end{align}
 with $M_\infty^{(2k+1,\theta)} $ given by
 \begin{align}\label{Martingale-General-case}
 M_\infty^{(2k+1,\theta)}:= \lim_{t\to\infty} e^{-(\beta(\mu-1)-\frac{\theta^2}{2})t}\sum_{u\in N(t)}  1_{\{\min_{s\leq t} X_u(s)> 0\}}e^{\theta X_u(t)} t^{(2k+1)/2} H_{2k+1}\left(\frac{X_u(t)}{\sqrt{t}}\right).
 \end{align}

The counterpart of Theorem \ref{thrm2} in the general case is as follows:
	For any given $m\in \N$, if \eqref{LlogL} holds for  some
		  $\lambda >6m +6$, then
	
	(i)  for any $x>0, b>a\geq 0 $, $\P_x$-almost surely, as $t\to\infty$,
	\begin{align}\label{step_44}
		&\frac{Z_{t}^{(0,\infty)} ((a,b))}{t^{-3/2}e^{\beta (\mu-1)t}}\nonumber\\
		&=- \sqrt{\frac{2}{\pi}}\sum_{\ell=0}^m \frac{H_{2\ell+2}(0)}{t^{\ell}} \sum_{k=0}^{\ell} \frac{M_\infty^{(2k+1,0)} }{(2k+1)! (2\ell -2k+1)!}
		\int_a^b z^{2\ell -2k+1}  \mathrm{d} z+o(t^{-m});
	\end{align}
	
	(ii) for any $x>0, a\geq 0 $, $\P_x$-almost surely, as $t\to\infty$,
	\begin{align}\label{step_45}
		&\frac{Z_{t}^{(0,\infty)} ((a,\infty))}{t^{-1/2}e^{\beta(\mu-1)t}}= \sqrt{\frac{2}{\pi}}\sum_{\ell=0}^m \frac{H_{2\ell}(0)}{t^{\ell}} \sum_{k=0}^{\ell} \frac{M_\infty^{(2k+1,0)} }{(2k+1)! (2\ell -2k)!} a^{2\ell-2k} +o(t^{-m}),
	\end{align}
 with $M_\infty^{(2k+1,0)} $ given in \eqref{Martingale-General-case}.
\end{remark}

Note that if $\sum^\infty_{k=1}k^{1+\epsilon}p_k<\infty$ for some $\epsilon>0$, then the conclusions of
Theorems \ref{thm1} and \ref{thrm2} hold for all $m\in \N$.

Our strategy for proving these two theorems is as follows. We will choose appropriate $\kappa >1$ and define
\begin{align}
	r_n:= n^{\frac{1}{\kappa}}, \quad n\in \mathbb{N}.
\end{align}
 We first study the asymptotic expansion along $\{r_n: n\in \N\}$,
which is given by Proposition \ref{Discrete-time-expansion}, and then control the behavior for $t\in (r_n, r_{n+1})$, see
 Lemma \ref{Continuous-time-approximation}
below. Once we have Proposition \ref{Discrete-time-expansion} and Lemma \ref{Continuous-time-approximation}, the proofs
of Theorems \ref{thm1} and  \ref{thrm2} are straight-forward.
To prove Proposition \ref{Discrete-time-expansion}, we first show that $Z_{r_n}^{(0,\infty)}(A) \approx \mathbb{E}_x\left( Z_{r_n}^{(0,\infty)}(A) \big| \mathcal{F}_{\sqrt{r_n}}\right)$, where $\mathcal{F}_t$ is the $\sigma$-field generated
by the branching Brownian motion with drift up to $t$, see Lemma \ref{lemma2} below for more details.
We then prove Proposition \ref{Discrete-time-expansion} with the help of two series expansions for the normal distribution (see Lemma \ref{lemma4} and Lemma \ref{lemma1}) and the convergence rate for martingales
$M_t^{(2k+1,\theta)}$
 (see Proposition \ref{Convergence-rate-martingale}).
To prove Lemma \ref{Continuous-time-approximation},
we first give a lower bound of $Z_t^{(0,\infty)}(A)$ (see Lemma \ref{Continuous-Lower-bound}), then accomplish the
proof of Lemma \ref{Continuous-time-approximation} by proving an upper bound of $Z_t^{(0,\infty)}((0,\infty))$.

\section{Preliminaries}\label{pre}

\subsection{Spine decomposition}

Define
\begin{align}\label{Change-of-measure}
	\frac{\mathrm{d}\mathbf{P}_x}{\mathrm{d} \mathbb{P}_x} \bigg|_{\mathcal{F}_t}
	:= \frac{\sum_{u\in N(t)} 1}{e^t},
\end{align}
then, under $\mathbf{P}_x$,
 the branching Brownian motion has the following spine decomposition(see \cite{Ky04}, or \cite{RS20} for a more general case):

(i) there is an initial marked particle at $x \in \mathbb{R}$ which moves according to standard Brownian motion with drift $-\theta$;

(ii) the branching rate of this marked particle is $2$;

(iii) when the marked particle dies at site $y$,
it gives birth to $\widehat{L}$ children with
$\mathbf{P}_x(\widehat{L}=k)=kp_k/2$;

(iv)  one of these children is uniformly selected and marked, and
the marked child evolves as its parent and the other children evolve  with law $\mathbb{P}_y$, where $\mathbb{P}_y$ denotes the law of a branching Brownian motion starting at $y$, and all the children evolve independently.

We use $\xi_t$  and $X_\xi(t)$ to denote the marked particle at time $t$ and the position of this marked particle respectively.
 By \cite[Theorem 2.11]{RS20}, we can get that for $u\in N(t)$,
\begin{align}\label{Many-to-one}
	\mathbf{P}_x
	\left(\xi_t =u \big| \mathcal{F}_t\right) = \frac{1}{\sum_{u\in N(t)} 1}.
\end{align}
Using \eqref{Many-to-one}, we get the following many-to-one formula:

\begin{lemma}\label{General-many-to-one}
		For any $x\in \R$, $t>0$ and $u\in N(t)$, let $\Gamma(u, t)$ be a  non-negative $\mathcal{F}_t$-measurable  random variable. Then
		\[
		\E_x \Big(\sum_{u\in N(t)} \Gamma(u,t) \Big) =e^t
		\mathbf{E}_x \left(\Gamma (\xi_t, t)\right).
		\]
\end{lemma}
\textbf{Proof: } Combining \eqref{Change-of-measure} and \eqref{Many-to-one}, we get
\begin{align}
&	\E_x \Big(\sum_{u\in N(t)} \Gamma(u,t) \Big)
	= e^{t}\mathbf{E}_x \Big(\sum_{u\in N(t)} \frac{\Gamma(u,t)}{\sum_{v\in N(t)} 1} \Big) = e^t \mathbf{E}_x \Big(\sum_{u\in N(t)}  \Gamma(u,t)  	 \mathbf{P}_x \Big(\xi_t =u \big| \mathcal{F}_t\Big)\Big)\nonumber\\
	& = e^t  \mathbf{E}_x\Big(     \mathbf{E}_x\Big(  \sum_{u\in N(t)} 1_{ \left\{\xi_t =u \right\}} \Gamma(u,t)  	   \big| \mathcal{F}_t\Big)\Big)
=  e^t \mathbf{E}_x\Big(    \Gamma (\xi_t,t)  	  \sum_{u\in N(t)} 1_{ \left\{\xi_t =u \right\}} \Big) = e^t\mathbf{E}_x\left(\Gamma(\xi_t, t)\right).
\end{align}
\hfill$\Box$

\subsection{Some useful facts}

\begin{lemma}\label{Useful-Ineq-1}
	(i) Let $\ell \in [1,2]$.
	Then for any
	finite family of independent centered random variables  $\left\{X_i : i=1, \dots, n\right\}$ with ${\rm E}|X_i|^\ell <\infty$ for all $i=1, \dots, n$,
	it holds that
	\[
	{\rm E} 	\big|\sum^n_{i=1} X_i\big|^\ell	\le 2 \sum^n_{i=1} {\rm E} |X_i|^\ell.
	\]
	(ii) For any $\ell\in [1,2]$ and any random variable $X$ with ${\rm E}|X|^\ell <\infty$,
	\[
	{\rm E}\left|X- {\rm E} X\right|^\ell \lesssim {\rm E} |X|^\ell \leq ({\rm E} X^2)^{\ell /2}.
	\]
\end{lemma}
\textbf{Proof: } For (i), see  \cite[Theorem 2]{BE65}. (ii) follows easily from Jensen's inequality.
\hfill$\Box$

We will use $\Pi_x^{-\theta}$ and $\Pi_x$ to denote the laws of a Brownian motion with drift $-\theta$ starting from $x$ and a standard Brownian motion starting from $x$ respectively.
Let $\phi(x):= \frac{1}{\sqrt{2\pi}}e^{-x^2/2}$
and $\Phi(x):=\int_{-\infty}^x \phi(y)\mathrm{d}y$.
\begin{lemma}\label{lemma3}
	(i) Let  $(B_t, \Pi_x^\uparrow)$ be a $3$-dimensional Bessel process with transition density
	\begin{align}\label{Transition-Probability}
		p_t^\uparrow(x,y)= 1_{\{y>0\}}\frac{ye^{-(y-x)^2/(2t)}}{x\sqrt{2\pi t}}\left(1-e^{-2xy/t}\right)=1_{\{y>0\}} \frac{y}{x\sqrt{t}}\left(\phi\left( \frac{y-x}{\sqrt{t}}\right) -\phi\left( \frac{y+x}{\sqrt{t}}\right)  \right).
	\end{align}
	Then for any $\theta \in [0,\sqrt{2})$, we have
	\begin{align}\label{step_13}
		\frac{\mathrm{d}\Pi_x^\uparrow}{\mathrm{d}\Pi_x^{-\theta}}\bigg|_{\sigma(B_s: s\leq t)}  =\frac{B_t e^{\theta(B_t-x)+\frac{\theta^2}{2}t}}{x}1_{\{\min_{s\leq t} B_s >0\}}.
	\end{align}
	(ii) For any $t,x,y>0$,
	\begin{align}\label{Upper-Bound-Transition-Prob}		
	p_t^\uparrow(x,y)
	\lesssim \frac{y^2}{t^{3/2}} e^{-(x-y)^2/(2t)}\leq \frac{y^2}{t^{3/2}}.
	\end{align}
\end{lemma}
\textbf{Proof: }
 (ii) follows from the  inequality $1-e^{-x}\leq x, x>0$. For (i),
note that under $\Pi_x$, both
\[
\frac{B_t}{x}1_{\{\min_{s\leq t} B_s >0\}} \quad \mbox{and}\quad e^{-\theta(B_t -x)-\frac{\theta^2}{2}t}
\]
are mean $1$ non-negative martingales and that we have the following change-of-measure:
\[
\frac{\mathrm{d}\Pi_x^\uparrow}{\mathrm{d}\Pi_x}\bigg|_{\sigma(B_s: s\leq t)}  = \frac{B_t}{x}1_{\{\min_{s\leq t} B_s >0\}}\quad \mbox{and} \quad \frac{\mathrm{d}\Pi_x^{-\theta}}{\mathrm{d} \Pi_x}\bigg|_{\sigma(B_s: s\leq t)}  = e^{-\theta(B_t -x)-\frac{\theta^2}{2}t}.
\]
Therefore, (i) follows from
\[
 \frac{\mathrm{d}\Pi_x^\uparrow}{\mathrm{d}\Pi_x^{-\theta}}\bigg|_{\sigma(B_s: s\leq t)}  = \frac{ \mathrm{d}\Pi_x^\uparrow}{\mathrm{d}\Pi_x}\bigg|_{\sigma(B_s: s\leq t)} \times \left( \frac{\mathrm{d}\Pi_x^{-\theta}}{ \mathrm{d}\Pi_x}\bigg|_{\sigma(B_s:s\leq t)} \right)^{-1} =\frac{B_t e^{\theta(B_t-x)+\frac{\theta^2}{2}t}}{x}1_{\{\min_{s\leq t} B_s >0\}}.
 \]
\hfill$\Box$

\begin{lemma}\label{Sharp-Probability}
	(i) Let $\theta \in (0,\sqrt{2})$. For any $x,t>0$ and 	Borel set $A\subset (0,\infty)$, it holds that
	\begin{align}
		\Pi_x^{-\theta}\left(\min_{s\leq t} B_s >0, B_t\in A\right)
		= \sqrt{\frac{2}{\pi}}
                 xe^{\theta x} t^{-3/2}e^{-\frac{\theta^2}{2}t}
                 \left(\int_A ye^{-\theta y}\mathrm{d} y+ \varepsilon_A(x,t) \right),
	\end{align}
 with $\varepsilon_A(x,t)$ satisfying
 \[
 \left|\varepsilon_A(x,t)\right|\leq C_\theta \left(1\land \frac{(x+1)^2}{t}\right)
 \]
 for some constant $C_\theta$ depending on $\theta$ only. In particular, for any fixed $\theta\in (0,\sqrt{2})$,
 \[
	\Pi_x^{-\theta}\Big(\min_{s\leq t} B_s >0, B_t\in A\Big)
	\lesssim xe^{\theta x}t^{-3/2}e^{-\frac{\theta^2}{2}t}.
 \]
(ii) For any Borel set
$A\subset (0,\infty)$ and $x,t>0$, it holds that
 \begin{align}
		& \Pi_x\Big(\min_{s\leq t} B_s >0, B_t\in A\Big)\lesssim  xt^{-1/2} 1_{\{\sup\{y: y\in A\}=\infty\}}  + x t^{-3/2} 1_{\{\sup\{y: y\in A\}<\infty\}}.
 \end{align}

 \noindent
 (iii) 	For any $x,t>0$, $\theta \in [0,\sqrt{2})$ and Borel set  $A\subset (0,\infty)$, it holds that
 \[
 \Pi_x^{-\theta}\left(B_t\in A\right)\leq e^{\theta x}e^{-\frac{\theta^2}{2}t}.
 \]
\end{lemma}
\textbf{Proof: } For (i), see \cite[Lemma 3.1]{LS20}; for (ii),
when $\sup\{y: y\in A\}=\infty$,
by the reflection principle for Brownian motion, we have
\[
		\Pi_x\Big(\min_{s\leq t} B_s >0\Big) = \Pi_0(|B_t|\leq x)
	= 2\int_0^x \frac{1}{\sqrt{2\pi t}}e^{-y^2/(2t)}\mathrm{d} y \lesssim
	 \frac{x}{\sqrt{t}}.
\]
 When $\sup\{y: y\in A\}<\infty$, by Lemma \ref{lemma3},
\begin{align}
		\Pi_x\Big(\min_{s\leq t} B_s >0, B_t\in A\Big)= x\Pi_x^\uparrow\left(\frac{1_{A}(B_t)}{B_t}\right)= x\int_A
		\frac{1}{y}p_t^\uparrow(x,y)\mathrm{d}y
		 \lesssim \frac{x}{t^{3/2}}\int_A y\mathrm{d}y\lesssim \frac{x}{t^{3/2}}.
\end{align}
For (iii),  by Girsanov's theorem,
\begin{align}
		\Pi_x^{-\theta}\left(B_t\in A\right)\leq 	
		\Pi_x^{-\theta}\left(B_t>0\right)= \Pi_x \Big(e^{-\theta(B_t-x)-\frac{\theta^2}{2}t}1_{\{B_t>0\}}\Big)\leq e^{\theta x} e^{-\frac{\theta^2}{2}t}.
\end{align}
\hfill$\Box$

\begin{lemma}\label{Properties-of-Hermite-Polynomial}
	(i) For any $k\geq 1$ and $x\in \mathbb{R}$,
	\begin{align}\label{Upper-H_k-2}
		\left|H_k(x)\right|\leq 2\sqrt{k!}e^{x^2/4}.
	\end{align}
   Consequently, it holds that
   \[
   \sup_{y\in\R} \left| \phi(y) H_k(y)\right| \leq \sqrt{\frac{2}{\pi}}\sqrt{k!} \sup_{y\in\R} e^{-y^2/4} = \sqrt{\frac{2}{\pi}}\sqrt{k!}.
   \]
   (ii) For any $k\in\N$,  there exists a constant $C(k)$ such that for all $x\in\R$,
   \[
   |H_{2k+1}(x)|\leq C(k) |x|\left(|x|^{2k}+1\right).
   \]
\end{lemma}
\textbf{Proof:} For (i), see \cite[(4.1) ]{GL21}; (ii) follows from the definition of $H_{2k+1}(x)$.
\hfill$\Box$

\begin{lemma}\label{Asymptotic-expansion}
	For any $\rho\in (0,1), b,x \in \mathbb{R}$, it holds that
	\[
	\Phi\left(\frac{b-\rho x}{\sqrt{1-\rho^2}}\right) = \Phi (b)- \phi(b) \sum_{k=1}^\infty \frac{\rho^k}{k!} H_{k-1}(b) H_k(x).
	\]
\end{lemma}
\textbf{Proof: } See \cite[Lemma 4.2.]{GL21}.
\hfill$\Box$

 Recall that $r_n=n^{1/\kappa}$.
 Applying Lemma \ref{Asymptotic-expansion} with $\rho= r_n^{-1/4}$, $b=r_n^{-1/2}z$ and
$x= r_n^{-1/4}y$,
we get that for any $z, y\in \R$,
\[
	\Phi\left(\frac{z-y}{\sqrt{r_n-\sqrt{r_n}}}\right) = \Phi \left(\frac{z}{\sqrt{r_n}}\right)- \phi\left(\frac{z}{\sqrt{r_n}}\right) \sum_{k=1}^\infty \frac{1}{k!} \frac{1}{r_n^{k/2}} H_{k-1}\left(\frac{z}{\sqrt{r_n}}\right) r_n^{k/4}H_k\left(\frac{y}{r_n^{1/4}}\right).
\]
Noting that, for any $k\in\N$,
$H_{2k}$
is an even function and $H_{2k+1}$ is an odd function, we get that
\begin{align}\label{Asymptotic-expansion-4}
	&\Phi\left(\frac{z+y}{\sqrt{r_n-\sqrt{r_n}}}\right) - \Phi\left(\frac{z-y}{\sqrt{r_n-\sqrt{r_n}}}\right) \nonumber\\
	&= 2\phi\left(\frac{z}{\sqrt{r_n}}\right)  \sum_{k=0}^\infty \frac{1}{(2k+1)!} \frac{1}{r_n^{(2k+1)/2}} H_{2k}\left(\frac{z}{\sqrt{r_n}}\right)  r_n^{(2k+1)/4}H_{2k+1}\left(\frac{y}{r_n^{1/4}}\right) .
\end{align}
\begin{lemma}\label{lemma4}
For any given $m\in \mathbb{N}$ and $\kappa >1$, let $K>0$ be a fixed constant and $J$ be an integer such that $J > 2m+ \frac{K\kappa -1}{2} $.
Then for any $y, z\in \R$, it holds that
	\begin{align}\label{e:neweqn}
			&\Phi\left(\frac{z+y}{\sqrt{r_n-\sqrt{r_n}}}\right) - \Phi\left(\frac{z-y}{\sqrt{r_n-\sqrt{r_n}}}\right) \nonumber\\
		&= 2\phi\left(\frac{z}{\sqrt{r_n}}\right)  \sum_{k=0}^J  \frac{1}{(2k+1)!} \frac{1}{r_n^{(2k+1)/2}} H_{2k}\left(\frac{z}{\sqrt{r_n}}\right)  r_n^{(2k+1)/4}H_{2k+1}\left(\frac{y}{r_n^{1/4}}\right) + \varepsilon_{m,y, z, n,\kappa}
	\end{align}
with
	\[
	r_n^{(2m+1)/2} \sup\left\{ \left|\varepsilon_{m,y,z, n,\kappa}\right|: \ z\in \R ,\ |y| \leq \sqrt{K\sqrt{r_n}\log n} \right\}
	\stackrel{n\to\infty}{\longrightarrow}0.
	\]
\end{lemma}
\textbf{Proof: } By Lemma \ref{Properties-of-Hermite-Polynomial} (i), for all $k\geq 1$, $z\in\R$ and $|y|\leq \sqrt{K\sqrt{r_n}\log n}$,
\begin{align}
	& \frac{r_n^{(2m+1)/2}}{(2k+1)!}  \frac{1}{r_n^{(2k+1)/2}} \left|\phi\left(\frac{z}{\sqrt{r_n}}\right) H_{2k}\left(\frac{z}{\sqrt{r_n}}\right) \right|  \times  \left|r_n^{(2k+1)/4} H_{2k+1}\left(\frac{y}{r_n^{1/4}}\right) \right|\\
	& \leq \frac{4 r_n^{(2m+1)/2}}{r_n^{(2k+1)/4}} \frac{1}{\sqrt{2\pi}}e^{y^2/\sqrt{r_n}}\leq \frac{4}{\sqrt{2\pi}} \frac{1}{n^{(2k-1-4m)/(4\kappa)}} n^{K/4}.
\end{align}
Combining this with \eqref{Asymptotic-expansion-4}, we get that \eqref{e:neweqn}  holds with
\begin{align}
	& r_n^{m} \sup_{z\in \R, |y|\leq \sqrt{K\sqrt{r_n}\log n} } \left|\varepsilon_{m,y,z, n}\right| \nonumber\\
	&\leq\frac{4}{\sqrt{2\pi}}  \sum_{k=J+1}^\infty \frac{1}{n^{(2k-1-4m)/(4\kappa)}} n^{K/4}\lesssim  \frac{1}{n^{(2(J+1)-1-4m)/(4\kappa)}} n^{K/4},
\end{align}
which tends to $0$ is since $J > 2m+ \frac{K\kappa -1}{2}$.
\hfill$\Box$

Taking derivative with respect to $b$ in Lemma \ref{Asymptotic-expansion},  and using the fact that
\begin{align}\label{K-th-derivative}
	\frac{\mathrm{d}^k}{\mathrm{d} b^k} \Phi(b)=(-1)^{k-1} H_{k-1}(b)\phi(b),
\end{align}
we get that

\begin{align}\label{Asymptotic-expansion-1}
	\frac{1}{\sqrt{1-\rho^2}} \phi \left( \frac{b-\rho x}{\sqrt{1-\rho^2}}\right) = \phi(b) + \phi(b) \sum_{k=1}^\infty \frac{\rho^k}{k!}H_k(b)H_k(x).
\end{align}
 Taking
 $\rho= r_n^{-1/4}$, $b=r_n^{-1/2}z$ and
$x= r_n^{-1/4}y$ in \eqref{Asymptotic-expansion-1}, we get that for any $z, y\in \R$,
\begin{align}\label{Asymptotic-expansion-2}
	\frac{\sqrt{r_n}}{\sqrt{r_n-\sqrt{r_n}}}\phi\left( \frac{z-y}{\sqrt{r_n-\sqrt{r_n}}}\right) = \phi\left(\frac{z}{\sqrt{r_n}}\right) \left(1+ \sum_{k=1}^\infty \frac{1}{k!} \frac{1}{r_n^{k/2}} H_k\left(\frac{z}{\sqrt{r_n}}\right)  r_n^{k/4}H_k\left(\frac{y}{r_n^{1/4}}\right) \right).
\end{align}
Noting that, for any $k\in\N$,
$H_{2k}$ is an even function and $H_{2k+1}$ is an odd function,
we deduce from \eqref{Asymptotic-expansion-2} that
\begin{align}\label{Asymptotic-expansion-3}
	&\frac{\sqrt{r_n}}{\sqrt{r_n-\sqrt{r_n}}}\left( \phi\left( \frac{z-y}{\sqrt{r_n-\sqrt{r_n}}}\right) - \phi\left(\frac{z+y}{\sqrt{r_n-\sqrt{r_n}}}\right) \right)\nonumber\\
	&= 2\phi\left(\frac{z}{\sqrt{r_n}}\right)  \sum_{k=0}^\infty \frac{1}{(2k+1)!} \frac{1}{r_n^{(2k+1)/2}} H_{2k+1}\left(\frac{z}{\sqrt{r_n}}\right)  r_n^{(2k+1)/4}H_{2k+1}\left(\frac{y}{r_n^{1/4}}\right) .
\end{align}
Using an argument similar to that leading to Lemma \ref{lemma4}, we also have the following lemma. We omit the proof.
\begin{lemma}\label{lemma1}
For any given $m\in \mathbb{N}$ and $\kappa >1$, let $K>0$ be a fixed constant and $J$ be an integer such that $J > 2m+ \frac{K\kappa +1}{2} $. Then for any $y, z\in \R$, it holds that
	\begin{align}
	&\frac{\sqrt{r_n}}{\sqrt{r_n-\sqrt{r_n}}}\left( \phi\left( \frac{z-y}{\sqrt{r_n-\sqrt{r_n}}}\right) - \phi\left(\frac{z+y}{\sqrt{r_n-\sqrt{r_n}}}\right) \right)\nonumber\\
	&= 2\phi\left(\frac{z}{\sqrt{r_n}}\right)  \sum_{k=0}^J \frac{1}{(2k+1)!} \frac{1}{r_n^{(2k+1)/2}} H_{2k+1}\left(\frac{z}{\sqrt{r_n}}\right)  r_n^{(2k+1)/4}H_{2k+1}\left(\frac{y}{r_n^{1/4}}\right) + \varepsilon_{m,y, z, n,\kappa}
	\end{align}
with
	\[
	r_n^{m+1} \sup\left\{ \left|\varepsilon_{m,y,z, n,\kappa}\right|: \ z\in \R ,\ |y| \leq \sqrt{K\sqrt{r_n}\log n} \right\}
	\stackrel{n\to\infty}{\longrightarrow}0.
	\]
\end{lemma}

\section{Proofs of the main results}\label{Main}

\subsection{Convergence rate for the martingales}

\begin{prop}
\label{Convergence-rate-martingale}
Suppose $x>0$ and $\theta\in[0,\sqrt{2})$
	(i)
	For any $k\in \N$,
$\{M_t^{(2k+1, \theta )}, t\geq 0;\, \mathbb{P}_x\}$ is a martingale.
	(ii) If \eqref{LlogL} holds for  some $\lambda > k$, then $M_t^{(2k+1,\theta)}$ converges to a limit $M_\infty^{(2k+1,\theta)}$ $\mathbb{P}_x$-a.s. and in $L^1(\P_x)$.
Moreover,  for any $\eta\in (0, \lambda - k)$,
	as $t\to\infty$,
	\[
	M_t^{(2k+1,\theta)} - M_\infty^{(2k+1,\theta)}=  o(t^{-(\lambda -k)+\eta}),\quad
\mathbb{P}_x\mbox{-a.s.}
	\]
\end{prop}
\textbf{Proof: } (i) We will use $v<u$ to denote that $v$ is an ancestor of $u$ and $v\leq u$ to denote $v=u$ or $v<u$. By the Markov property and Lemma \ref{General-many-to-one}, for any $t,s>0$,
\begin{align}\label{step_20}
	&\mathbb{E}_x\left(M_{s+t}^{(2k+1,\theta)}\big| \mathcal{F}_t\right)=e^{-(1-\frac{\theta^2}{2})(t+s)}\sum_{v\in N(t)}  1_{\{\min_{r\leq t} X_v(r)> 0\}}(t+s)^{(2k+1)/2} \nonumber\\
	&\quad\quad\times  \mathbb{E}_{X_v(t)}\left(\sum_{u\in N(t+s): v\leq u} 1_{\{\min_{r\leq s} X_v(t+r) > 0\}}e^{\theta X_v(t+s)}  H_{2k+1}\left(\frac{X_u(t+s)}{\sqrt{t+s}} \right) \Big| \mathcal{F}_t\right)\nonumber\\
	& =e^{-(1-\frac{\theta^2}{2})(t+s)}\sum_{v\in N(t)}  1_{\{\min_{r\leq t} X_v(r)> 0\}} \nonumber\\
	&\quad\quad\times e^{s}
	\mathbf{E}_{X_v(t)}
	\left(1_{\{\min_{r\leq s} X_\xi(r) > 0\}}e^{\theta X_\xi (s)}  (t+s)^{(2k+1)/2}H_{2k+1}\left(\frac{X_\xi(s)}{\sqrt{t+s}} \right) \right)
	\nonumber \\	& =: e^{-(1-\frac{\theta^2}{2})(t+s)}\sum_{v\in N(t)}  1_{\{\min_{r\leq t} X_v(r)> 0\}} F(s, t, X_v(t)).
\end{align}
Note that
$X_\xi(s)$ under  $\mathbf{P}_x$  is a standard Brownian motion with drift $-\theta$. It follows from Lemma \ref{lemma3} (i) that
\begin{align}
	& F(s, t, X_v(t)) = e^{(1-\frac{\theta^2}{2})s } X_v(t)e^{\theta X_v(t)}
	 \Pi_{X_v(t)}^\uparrow
	\left( \frac{(t+s)^{(2k+1)/2}}{B_s}H_{2k+1}\left(\frac{B_s}{\sqrt{t+s}} \right) \right)\nonumber\\
	& = e^{(1-\frac{\theta^2}{2})s } X_v(t)e^{\theta X_v(t)} \int_0^\infty \frac{(t+s)^{(2k+1)/2}}{y}H_{2k+1}\left(\frac{y}{\sqrt{t+s}}\right)\nonumber\\
	&\quad\quad \times \frac{y}{X_v(t)\sqrt{s}}\left(\phi\left( \frac{y-X_v(t)}{\sqrt{s}}\right) -\phi\left( \frac{y+X_v(t)}{\sqrt{s}}\right) \right) \mathrm{d}y\nonumber\\
	& = e^{(1-\frac{\theta^2}{2})s } e^{\theta X_v(t)} \int_0^\infty (t+s)^{(2k+1)/2}H_{2k+1}\left(\frac{y}{\sqrt{t+s}}\right) \frac{\phi\left( \frac{y-X_v(t)}{\sqrt{s}}\right) -\phi\left( \frac{y+X_v(t)}{\sqrt{s}}\right) }{\sqrt{ s}}\mathrm{d}y.
\end{align}
Using the fact that
$H_{2k+1}(\cdot)$ is an odd function and that $\phi(\cdot)$ is an even function,
we have
\begin{align}
	&  F(s, t, X_v(t)) = e^{(1-\frac{\theta^2}{2})s } e^{\theta X_v(t)} \int_{-\infty}^\infty
	 (t+s)^{(2k+1)/2}H_{2k+1}\left(\frac{y}{\sqrt{t+s}}\right) \frac{\phi\left( \frac{y-X_v(t)}{\sqrt{s}}\right) }{\sqrt{ s}}\mathrm{d}y\nonumber\\
	& = e^{(1-\frac{\theta^2}{2})s } e^{\theta X_v(t)} \int_{-\infty}^\infty (t+s)^{(2k+1)/2}H_{2k+1}\left(z\right) \frac{\sqrt{t+s}}{\sqrt{s}} \phi \left( \frac{ \sqrt{t+s} z-X_v(t)}{\sqrt{s}}\right) \mathrm{d}z.
\end{align}
Taking $\rho:= \sqrt{t/(s+t)}, b=z, x= X_v(t)/\sqrt{t}$ in \eqref{Asymptotic-expansion-1}, we see that
\[
	\frac{\sqrt{t+s}}{\sqrt{s}} \phi \left( \frac{ \sqrt{t+s} z-X_v(t)}{\sqrt{s}}\right) = \phi(z) + \phi(z) \sum_{\ell =1}^\infty \frac{t^{\ell /2}}{\ell !}\frac{1}{(t+s)^{\ell /2}}H_\ell (z)H_\ell \left(\frac{X_v(t)}{\sqrt{t}}\right).
\]
Combining these with Lemma \ref{Properties-of-Hermite-Polynomial}(i) and the fact $\phi(z) e^{z^2/4}\lesssim e^{-z^2/4}$, we
can easily get that for any $s,t>0, y\in \R$, the series
\[
H_{2k+1}\left(z\right) \phi(z) \sum_{\ell =1}^\infty \frac{t^{\ell /2}}{\ell !}\frac{1}{(t+s)^{\ell /2}}H_\ell (z)H_\ell \left(\frac{y}{\sqrt{t}}\right)
\]
is uniformly convergent in $z\in \R$. Now applying the property $\int_{-\infty}^\infty H_{m}(z)H_n(z)\phi(z)\mathrm{d}z= \delta_{m,n} n!$,
we get
\begin{align}
	&F(s, t, X_v(t)) \nonumber\\
	&=e^{(1-\frac{\theta^2}{2})s } e^{\theta X_v(t)} \int_{-\infty}^\infty (t+s)^{(2k+1)/2}\frac{t^{(2k+1) /2}}{(2k+1)!}\frac{H_{2k+1}\left(z\right) \phi(z)}{(t+s)^{(2k+1) /2}}H_{2k+1} (z)H_{2k+1} \left(\frac{X_v(t)}{\sqrt{t}}\right) \mathrm{d}z\nonumber\\
	& = e^{(1-\frac{\theta^2}{2})s } e^{\theta X_v(t)} t^{(2k+1) /2}H_{2k+1} \left(\frac{X_v(t)}{\sqrt{t}}\right).
\end{align}
Plugging this into \eqref{step_20}, we get (i).

(ii) Suppose \eqref{LlogL} holds for  some $\lambda > k$. If the first assertion of (ii) holds along $t=n\in \N$, then it is valid along all $t$
since for $t\in (n,n+1)$, $M_t^{(2k+1, \theta)}=\E_x\left(M_{n+1}^{(2k+1,\theta)}\big| \mathcal{F}_t\right)= \E_x\left(M_\infty^{(2k+1,\theta)}\big| \mathcal{F}_t\right)$.
In the following we  use two steps to prove the assertion of (ii) holds for the case $t=n\in \N$.

{\bf Step 1:}
In this step, we will define a truncated process $M_{n+1}^{(2k+1,\theta, B)}$ and give first moment estimate for $M_{n+1}^{(2k+1,\theta)}- M_{n+1}^{(2k+1,\theta, B)}$.
Let $d_v, O_v$ denote the death time and the number of offspring of $v$ respectively.  For $u\in N(n+1)$, let
$B_{n,u}$ be the event that, for all $v<u$ with $d_v\in (n, n+1)$, it holds that $O_v \leq e^{c_0n}$,
where $c_0>0$ is a small constant to be determined later. Define
\begin{align}
		&M_{n+1}^{(2k+1,\theta, B)}\nonumber\\
		&:= e^{-(1-\frac{\theta^2}{2})(n+1)}\sum_{u\in N(n+1)}  1_{\{\min_{s\leq n+1} X_u(s)> 0\}}e^{\theta X_u(n+1)} (n+1)^{(2k+1)/2} H_{2k+1}\left(\frac{X_u(n+1)}{\sqrt{n+1}}\right)1_{B_{n,u}}.
\end{align}
By the branching property, it holds that
\begin{align}
	& \left|M_{n+1}^{(2k+1,\theta)}- M_{n+1}^{(2k+1,\theta, B)}\right|\nonumber\\
	&\leq  e^{-(1-\frac{\theta^2}{2})(n+1)}\sum_{u\in N(n+1)}  1_{\{\min_{s\leq n+1} X_u(s)> 0\}}e^{\theta X_u(n+1)} (n+1)^{(2k+1)/2} \left| H_{2k+1}\left(\frac{X_u(n+1)}{\sqrt{n+1}}\right) \right| 1_{B_{n,u}^c}\nonumber\\
	& = e^{-(1-\frac{\theta^2}{2})(n+1)} \sum_{v\in N(n)}  1_{\{\min_{s\leq n} X_v(s)> 0\}}\nonumber\\
	&\quad\times \sum_{u\in N(n+1): v\leq u} 1_{\{\min_{s\leq 1} X_u(n+s)>0\}} e^{\theta X_u(n+1)}  (n+1)^{(2k+1)/2} \left| H_{2k+1}\left(\frac{X_u(n+1)}{\sqrt{n+1}}\right) \right| 1_{B_{n,u}^c}.
\end{align}
Using the Markov property first and then Lemma \ref{Properties-of-Hermite-Polynomial}(ii), we get that
\begin{align}
	&\mathbb{E}_x\left( \left|M_{n+1}^{(2k+1,\theta)}- M_{n+1}^{(2k+1,\theta, B)}\right| \Big|\mathcal{F}_n\right)\nonumber\\
	& \leq e^{-(1-\frac{\theta^2}{2})(n+1)} \sum_{v\in N(n)}  1_{\{\min_{s\leq n} X_v(s)> 0\}} \nonumber\\
	& \quad \times \mathbb{E}_{X_v(n)}\left(\sum_{u\in N(1)} 1_{\{\min_{s\leq 1} X_u(s)>0\}} e^{\theta X_u(1)}  (n+1)^{(2k+1)/2} \left| H_{2k+1}\left(\frac{X_u(1)}{\sqrt{n+1}}\right) \right| 1_{D_{n,u}^c} \right)\nonumber\\
	& \lesssim e^{-(1-\frac{\theta^2}{2})(n+1)} \sum_{v\in N(n)}  1_{\{\min_{s\leq n} X_v(s)> 0\}} \nonumber\\
	& \quad \times \mathbb{E}_{X_v(n)}\left(\sum_{u\in N(1)} 1_{\{\min_{s\leq 1} X_u(s)>0\}} X_u(1) e^{\theta X_u(1)}   \left((X_u(1))^{2k}+ (n+1)^k \right) 1_{D_{n,u}^c} \right),
\end{align}
where for $u\in N(1)$, $D_{n,u}$ denotes the event that,
for all $w<u$ with $d_w<1$,
it holds that $O_w\leq e^{c_0 n}$.
Let $d_i$ be the $i$-th splitting time of the spine and $O_i$ be the number of children produced by the spine at time $d_i$.
Define $D_{n, \xi_1}$ to be the event that, for all $i$ with $d_i<1$, it holds that $O_i\leq e^{c_0 n}$.
By Lemma \ref{General-many-to-one},
\begin{align}\label{step_21}
	&\mathbb{E}_x\left( \left|M_{n+1}^{(2k+1,\theta)}- M_{n+1}^{(2k+1,\theta, B)}\right| \Big|\mathcal{F}_n\right)\nonumber\\
	& \lesssim e^{-(1-\frac{\theta^2}{2})(n+1)} \sum_{v\in N(n)}  1_{\{\min_{s\leq n} X_v(s)> 0\}}\nonumber\\
	&\quad \times e \mathbf{E}_{X_v(n)}
	\left( 1_{\{\min_{s\leq 1} X_\xi (s)>0\}} X_\xi(1) e^{\theta X_\xi(1)}   \left(X_\xi(1)^{2k}+ (n+1)^k \right) 1_{D_{n,\xi_1}^c} \right)
	\nonumber\\ 	&=: e^{1-(1-\frac{\theta^2}{2})(n+1)} \sum_{v\in N(n)}  1_{\{\min_{s\leq n} X_v(s)> 0\}} R(X_v(n)).
\end{align}
Conditioned on the motion $X_\xi$, $\{d_i: i\geq 1\}$ are the atoms of a Poisson poin process with rate $2$ and $\{O_i: i\in \N\}$ are iid copies of $\widehat{L}$ with law
$\mathbf{P}_x(\widehat{L}=\ell)=\ell p_\ell/2$
 which are independent of $\{d_i: i\in \N\}$. Therefore, $D_{n,\xi_1}$ is independent of $X_\xi(t)$. Together with Lemma \ref{lemma3} (i), we get that
\begin{align}\label{step_22}
	&R(X_v(n))=  \mathbf{E}_{X_v(n)}
	\left( 1_{\{\min_{s\leq 1} X_\xi (s)>0\}} X_\xi(1) e^{\theta X_\xi(1)}   \left(X_\xi(1)^{2k}+ (n+1)^k \right) \right)
	\mathbf{E}_{X_v(n)}
	\left( 1_{D_{n,\xi_1}^c} \right)\nonumber\\
	& \leq  e^{-\frac{\theta^2}{2}} X_v(n)e^{\theta X_v(n)}\Pi_{X_v(n)}^\uparrow\left( B_1^{2k}+ (n+1)^k\right)
	\mathbf{E}_{X_v(n)}
	\left( \sum_{i:d_i\leq 1} 1_{\{O_i >e^{c_0 n}\}}\right)\nonumber\\
	&= 2e^{-\frac{\theta^2}{2}} X_v(n)e^{\theta X_v(n)}\Pi_{X_v(n)}^\uparrow\left( B_1^{2k}+ (n+1)^k\right)
	\mathbf{P}_x(\widehat{L}> e^{c_0n}).
\end{align}
Noticing that $(B_t, \Pi_x^\uparrow)$ is a $3$-dimensional Bessel process, we easily see that
\begin{align}\label{step_22'}
 \Pi_{X_v(n)}^\uparrow\left( B_1^{2k}+ (n+1)^k\right)\leq \Pi_0^\uparrow\left( (B_1+y)^{2k}+ (n+1)^k\right)\big|_{y= X_v(n)}\lesssim (X_v(n))^{2k}+n^k.
\end{align}
Noting that \eqref{LlogL} implies
$\mathbf{E}_x (\log_+^{1+\lambda} \widehat{L})<\infty$,
using \eqref{step_22} and \eqref{step_22'}, we obtain that
\begin{align}
	R(X_v(n))
	& \lesssim \frac{1}{n^{1+\lambda}}e^{-\frac{\theta^2}{2}} X_v(n)e^{\theta X_v(n)} \left( (X_v(n))^{2k}+n^k\right).
\end{align}
Plugging this inequality into \eqref{step_21}, we conclude that
\begin{align}
		&\mathbb{E}_x\left( \left|M_{n+1}^{(2k+1,\theta)}- M_{n+1}^{(2k+1,\theta, B)}\right| \Big|\mathcal{F}_n\right)\nonumber\\
		&\lesssim
		 \frac{e^{-(1-\frac{\theta^2}{2})n} }{n^{1+\lambda}}\sum_{v\in N(n)}  1_{\{\min_{s\leq n} X_v(s)> 0\}} X_v(n)e^{\theta X_v(n)} \left( (X_v(n))^{2k}+n^k\right).
\end{align}
Taking expectation with respect to $\P_x$,
applying Lemma \ref{General-many-to-one} first and then Lemma \ref{lemma3} (i), we get that
\begin{align}\label{First-moment-martingale}
	&\mathbb{E}_x\left( \left|M_{n+1}^{(2k+1,\theta)}- M_{n+1}^{(2k+1,\theta, B)}\right| \right)\nonumber\\
	& \lesssim \frac{e^{-(1-\frac{\theta^2}{2})n} }{n^{1+\lambda}} \mathbb{E}_x\Big(
	\sum_{v\in N(n)}  1_{\{\min_{s\leq n} X_v(s)> 0\}} X_v(n)e^{\theta X_v(n)} \left( (X_v(n))^{2k}+n^k\right)\Big)\nonumber\\
	& =
	\frac{e^{\frac{\theta^2}{2}n }}{n^{1+\lambda}}	\mathbf{E}_x
	\left( 1_{\{\min_{s\leq n} X_\xi(s)> 0\}}e^{\theta X_\xi(n)} X_\xi(n) \left( (X_\xi(n))^{2k}+n^k\right) \right)\nonumber\\
	& = \frac{xe^{\theta x}}{n^{1+\lambda}}\Pi_x^\uparrow \left((B_n)^{2k}
	+n^k\right) \leq  \frac{x e^{\theta x}}{n^{1+\lambda}}\Pi_0^\uparrow \left(	(B_n+x)^{2k}+n^k\right) \nonumber\\	&\lesssim \frac{1}{n^{1+\lambda}}\Pi_0^\uparrow \left((B_n)^{2k}
	+n^k\right) =\frac{n^k}{n^{1+\lambda}}\Pi_0^\uparrow \left((B_1)^{2k}	+1\right),
\end{align}
where in the last equality,
we used the fact that
$(B_t, \Pi_0^\uparrow)\stackrel{\mathrm{d}}{=} (\sqrt{t}B_1, \Pi_0^\uparrow).$

{\bf Step 2:}
In this step, we
will give an upper bound for the $\ell$-th moment of $M_{n+1}^{(2k+1,\theta, B)}- \mathbb{E}_x\left(M_{n+1}^{(2k+1,\theta, B)}\big|\mathcal{F}_n\right)$ for appropriate $\ell\in (1,2)$.
Combining this with Step 1 will yield the first result of (ii).
 For $v\in N(n)$, set
\[
J_{n,v}:= \sum_{u\in N(n+1): v\leq u} 1_{\{\min_{s\leq 1} X_u(n+s)>0\}} e^{\theta X_u(n+1)}  (n+1)^{(2k+1)/2} H_{2k+1}\left(\frac{X_u(n+1)}{\sqrt{n+1}}\right) 1_{B_{n,u}}.
\]
By the branching property, $\{J_{n,v}: v\in N(n)\}$ are independent conditioned on $\mathcal{F}_n$. Thus, for  any fixed $1< \ell < \min\{2/\theta^2, 2\}$ with $(\ell-1)^2\theta^2/ 2 < (\ell-1)(1-\frac{\theta^2}{2})$, by Lemma \ref{Useful-Ineq-1},
\begin{align}\label{step_27}
	& \mathbb{E}_x\left(\left|M_{n+1}^{(2k+1,\theta, B)}- \mathbb{E}_x\left(M_{n+1}^{(2k+1,\theta, B)}\big|\mathcal{F}_n\right)\right|^\ell \big|\mathcal{F}_n \right)\nonumber\\
	& = e^{-\ell (1-\frac{\theta^2}{2})(n+1)}\mathbb{E}_x\left(\left| \sum_{v\in N(n)}  1_{\{\min_{s\leq n} X_v(s)> 0\}}\left(J_{n,v}-\mathbb{E}_x\left(J_{n,v}\big|\mathcal{F}_n\right)\right)\right|^\ell \big|\mathcal{F}_n \right)\nonumber\\
	& \leq 2 e^{-\ell (1-\frac{\theta^2}{2})(n+1)} \sum_{v\in N(n)}  1_{\{ X_v(n)> 0\}} \mathbb{E}_x\left(\left|J_{n,v}-\mathbb{E}_x\left(J_{n,v}\big|\mathcal{F}_n\right)\right|^\ell \big|\mathcal{F}_n \right)\nonumber\\
	& \lesssim  e^{-\ell(1-\frac{\theta^2}{2})(n+1)} \sum_{v\in N(n)}  1_{\{ X_v(n)> 0\}} \left(\mathbb{E}_x\left( (J_{n,v})^2 \big|\mathcal{F}_n \right)\right)^{\ell /2}.
\end{align}
Define
\begin{align}
J_{n}^*:= \sum_{u\in N(1)} 1_{\{\min_{s\leq 1} X_u(s)>0\}} e^{\theta X_u(1)}  (n+1)^{(2k+1)/2}  \left| H_{2k+1}\left(\frac{X_u(1)}{\sqrt{n+1}}\right) \right| 1_{D_{n,u}}.
\end{align}
By the Markov property and Lemma \ref{General-many-to-one},
\begin{align}\label{step_26}
	&\mathbb{E}_x\left( (J_{n,v})^2 \big|\mathcal{F}_n \right)\leq \E_{X_v(n)}\left((J_n^*)^2\right)\nonumber\\
	&= \E_{X_v(n)}\left(  \sum_{u\in N(1)} 1_{\{\min_{s\leq 1} X_u(s)>0\}} e^{\theta X_u(1)}  (n+1)^{(2k+1)/2} \left| H_{2k+1}\left(\frac{X_u(1)}{\sqrt{n+1}}\right) \right| 1_{D_{n,u}}J_n^* \right)  \nonumber\\
	& = e 	\mathbf{E}_{X_v(n)}
	\left(   1_{\{\min_{s\leq 1} X_\xi(s)>0\}} e^{\theta X_\xi(1)}  (n+1)^{(2k+1)/2} \left| H_{2k+1}\left(\frac{X_\xi(1)}{\sqrt{n+1}}\right) \right| 1_{D_{n,\xi_1}}J_n^* \right).
\end{align}
Conditioned on $\{X_\xi, d_i, O_i: i\geq 1\},$ by the  Markov property,
on the event $D_{n,\xi_1}$, we have
\begin{align}\label{step_23}
	& \mathbf{E}_{X_v(n)}
	\left( J_n^* \big| X_\xi, d_i, O_i: i\geq 1\right)
	= \sum_{i: d_i\leq 1} 1_{\{\min_{s\leq d_i} X_\xi(s)>0\}} (O_i-1)\nonumber\\
	& \quad\quad\times \mathbb{E}_{X_\xi(d_i)}\left( \sum_{u\in N(z)} 1_{\{\min_{s\leq z} X_u(s)>0\}} e^{\theta X_u(z)}  (n+1)^{(2k+1)/2}  \left| H_{2k+1}\left(\frac{X_u(z)}{\sqrt{n+1}}\right) \right| 1_{D_{n,u}} \right)\big|_{z=1-d_i}\nonumber\\ &
	 \leq  \sum_{i: d_i\leq 1} 1_{\{\min_{s\leq d_i} X_\xi(s)>0\}} (O_i-1)\nonumber\\
	& \quad\quad\times \mathbb{E}_{X_\xi(d_i)}\left( \sum_{u\in N(z)} 1_{\{\min_{s\leq z} X_u(s)>0\}} e^{\theta X_u(z)}  (n+1)^{(2k+1)/2}  \left| H_{2k+1}\left(\frac{X_u(z)}{\sqrt{n+1}}\right) \right|  \right)\big|_{z=1-d_i}\nonumber\\
	& \leq e^{c_0n}  \sum_{i: d_i\leq 1} 1_{\{\min_{s\leq d_i} X_\xi(s)>0\}} e^{1-d_i} \nonumber\\
	&\quad\quad\times  	\mathbf{E}_{X_\xi(d_i)}
	\left( 1_{\{\min_{s\leq z} X_\xi(s)>0\}} e^{\theta X_\xi(z)}  (n+1)^{(2k+1)/2}  \left| H_{2k+1}\left(\frac{X_\xi(z)}{\sqrt{n+1}}\right) \right|  \right)\big|_{z=1-d_i}
	\nonumber\\ 	& =: e^{c_0n}  \sum_{i: d_i\leq 1} 1_{\{\min_{s\leq d_i} X_\xi(s)>0\}} e^{1-d_i}
	  F_n(X_\xi(d_i), 1-d_i),
\end{align}
where in the last inequality we used the fact that
$O_i-1\leq e^{c_0 n}$ on $D_{n,\xi_1}$ and Lemma \ref{General-many-to-one}.
Using Lemma \ref{Properties-of-Hermite-Polynomial}(ii) in the first inequality,
Lemma \ref{lemma3} (i) in the first equality
and an argument similar to that leading to  \eqref{step_22'} in the second inequality, we get
\begin{align}\label{step_24}
	&F_n(X_\xi(d_i), 1-d_i) \lesssim 	\mathbf{E}_{X_\xi(d_i)}
	\left( 1_{\{\min_{s\leq z} X_\xi(s)>0\}} e^{\theta X_\xi(z)}  X_\xi(z) \left((X_\xi(z))^{2k}+ (n+1)^k\right) \right)\big|_{z=1-d_i}\nonumber\\
	& = X_\xi(d_i)e^{\theta X_\xi(d_i)} e^{-\frac{\theta^2}{2}(1-d_i)} \Pi_{X_\xi(d_i)}^\uparrow \left( (B_z)^{2k}+ (n+1)^k\right)\big|_{z=1-d_i}\nonumber\\
	&\lesssim X_\xi(d_i)e^{\theta X_\xi(d_i)} \left( (X_\xi(d_i))^{2k}+n^k\right)
	 \sup_{z\in(0,1)}\Pi_0^\uparrow ((B_z)^{2k}+1)\nonumber\\
	&\lesssim  X_\xi(d_i)e^{\theta X_\xi(d_i)} \left( (X_\xi(d_i))^{2k}+n^k\right).
\end{align}
Combining \eqref{step_23} and \eqref{step_24}, we obtain that
\begin{align}\label{step_25}
	& \mathbf{E}_{X_v(n)} \left( J_n^* \big| X_\xi, d_i, O_i: i\geq 1\right)\nonumber\\
	& \lesssim e^{c_0 n} \sum_{i: d_i\leq 1} 1_{\{\min_{s\leq d_i} X_\xi(s)>0\}} e^{1-d_i}X_\xi(d_i)e^{\theta X_\xi(d_i)} \left( (X_\xi(d_i))^{2k}+n^k\right) \nonumber\\
	& \lesssim e^{c_0 n} \sum_{i: d_i\leq 1} 1_{\{\min_{s\leq d_i} X_\xi(s)>0\}} X_\xi(d_i)e^{\theta X_\xi(d_i)} \left( (X_\xi(d_i))^{2k}+n^k\right).
\end{align}
Plugging \eqref{step_25} into \eqref{step_26},
\begin{align}\label{step_25'}
	&\mathbb{E}_x\left( (J_{n,v})^2 \big|\mathcal{F}_n \right)\lesssim  e^{c_0n}
	\mathbf{E}_{X_v(n)}
	\bigg(   1_{\{\min_{s\leq 1} X_\xi(s)>0\}} e^{\theta X_\xi(1)}  (n+1)^{(2k+1)/2} \left| H_{2k+1}\left(\frac{X_\xi(1)}{\sqrt{n+1}}\right) \right|\nonumber\\
	&\quad\quad\times  \int_0^1 1_{\{\min_{s\leq r} X_\xi(s)>0\}} X_\xi(r)e^{\theta X_\xi(r)} \left( (X_\xi(r))^{2k}+n^k\right)\mathrm{d} r  \bigg)
	\nonumber \\ & \lesssim e^{c_0 n}
	\mathbf{E}_{X_v(n)}
	\bigg(   1_{\{\min_{s\leq 1} X_\xi(s)>0\}} e^{\theta X_\xi(1)}  X_\xi(1) \left((X_\xi(1))^{2k}+ (n+1)^k\right) \nonumber\\
	&\quad\quad\times  \int_0^1 1_{\{\min_{s\leq r} X_\xi(s)>0\}} X_\xi(r)e^{\theta X_\xi(r)} \left( (X_\xi(r))^{2k}+n^k\right)\mathrm{d} r  \bigg),
\end{align}
where  the last inequality follows by Lemma \ref{Properties-of-Hermite-Polynomial}(ii).
Using Lemma \ref{lemma3}, we can continue the estimate \eqref{step_25'} and get
\begin{align}\label{step_28}
	&\mathbb{E}_x\left( (J_{n,v})^2 \big|\mathcal{F}_n \right)
	\lesssim
	 X_v(n)e^{\theta X_v(n)} e^{c_0 n}\Pi_{X_v(n)}^\uparrow \bigg(   \left((B_1)^{2k}+ (n+1)^k\right) \int_0^1 B_r e^{\theta B_r} \left( (B_r)^{2k}+n^k\right)\mathrm{d} r  \bigg)\nonumber\\
	& \leq X_v(n)e^{2\theta X_v(n)} e^{c_0 n}\Pi_{0}^\uparrow \bigg(   \left((B_1 +y)^{2k}+ (n+1)^k\right) \int_0^1 (B_r+y) e^{\theta B_r} \left( (B_r+y)^{2k}+n^k\right)\mathrm{d} r  \bigg)\big|_{y=X_v(n)}\nonumber\\
	& \lesssim e^{c_0 n}(X_v(n)+1)^2e^{2\theta X_v(n)}\left((X_v(n))^{2k}+n^k\right)^2.
\end{align}
Combining   \eqref{step_27} and  \eqref{step_28}, we have
\begin{align}
	& \mathbb{E}_x\left(\left|M_{n+1}^{(2k+1,\theta, B)}- \mathbb{E}_x\left(M_{n+1}^{(2k+1,\theta, B)}\big|\mathcal{F}_n\right)\right|^\ell \big|\mathcal{F}_n \right)\nonumber\\
	& \lesssim e^{-\ell(1-\frac{\theta^2}{2})(n+1)} e^{c_0 \ell n/2} \sum_{v\in N(n)} 1_{\{ X_v(n)> 0\}} (X_v(n)+1)^{\ell }e^{ \ell \theta X_v(n)}\left((X_v(n))^{2k}+n^k\right)^\ell \nonumber\\
	&\leq e^{-\ell(1-\frac{\theta^2}{2})n} e^{c_0 \ell n/2} \sum_{v\in N(n)}  (|X_v(n)|+1)^{\ell }e^{ \ell \theta X_v(n)}\left((X_v(n))^{2k}+n^k\right)^\ell .
\end{align}
Taking expectation with respect to $\P_x$ and applying  Lemma \ref{General-many-to-one},
we conclude that
\begin{align}\label{Second-moment-martingale}
	& \mathbb{E}_x\left(\left|M_{n+1}^{(2k+1,\theta, B)}- \mathbb{E}_x\left(M_{n+1}^{(2k+1,\theta, B)}\big|\mathcal{F}_n\right)\right|^\ell  \right)\nonumber\\
	& \lesssim e^{-\ell(1-\frac{\theta^2}{2})n} e^{c_0 \ell n/2}  e^n
	\mathbf{E}_x
	 \left( (|X_\xi(n)|+1)^{\ell }e^{ \ell \theta X_\xi(n)}\left((X_\xi(n))^{2k}+n^k\right)^\ell \right)\nonumber\\
	& = e^{-\ell(1-\frac{\theta^2}{2})n} e^{c_0 \ell n/2}  e^n \Pi_x^{-\theta}\left((|B_n|+1)^{\ell }e^{ \ell \theta B_n}\left((B_n)^{2k}+n^k\right)^\ell \right)\nonumber\\
	& =e^{\theta x} e^{-(\ell-1)(1-\frac{\theta^2}{2})n} e^{c_0 \ell n/2}   \Pi_x\left((|B_n|+1)^{\ell }e^{ (\ell-1) \theta B_n}\left((B_n)^{2k}+n^k\right)^\ell \right)\nonumber\\
	& = e^{\theta \ell x}e^{-(\ell-1)(1-\frac{\theta^2}{2})n} e^{c_0 \ell n/2} e^{\frac{(\ell-1)^2\theta^2}{2}n}\Pi_x^{(\ell-1)\theta}\left((|B_n|+1)^{\ell }\left((B_n)^{2k}+n^k\right)^\ell \right)\nonumber\\
	&\lesssim e^{-(\ell-1)(1-\frac{\theta^2}{2})n} e^{c_0 \ell n/2} e^{\frac{(\ell-1)^2\theta^2}{2}n}
		    n^{2k\ell +\ell}=: n^{2k\ell +\ell}e^{-c_1 \ell n},
\end{align}
where in the second and third equalities
we used the change-of-measure $\frac{\mathrm{d}\Pi_x^{\eta}}{\mathrm{d}\Pi_x} \big|_{\sigma(B_s: s\leq t)}= e^{\eta (B_t-x)-\frac{\eta^2}{2}t}$ for $\eta= -\theta$ and $\eta =(\ell-1)\theta$ respectively.
Let $c_0>0$ be sufficiently small  so that $c_0\ell/2 <  (\ell-1)(1-\frac{\theta^2}{2})-(\ell-1)^2\theta^2/ 2 $, which implies that $c_1>0$.
Thus, using the inequality:
\begin{align}\label{New-step_2}
	& {\rm E}\left(|X- {\rm E} (X| \mathcal{F}) |\right) \leq {\rm E} \left(|X-Y|\right)+ {\rm E}\left( |Y- {\rm E} (Y| \mathcal{F}) |\right) + {\rm E} \left(\left|{\rm E}\left(X-Y \big| \mathcal{F}\right) \right| \right)\nonumber\\
	& \leq 2 {\rm E} \left(|X-Y|\right) + {\rm E}\left( |Y- {\rm E} (Y| \mathcal{F}) |^\ell \right) ^{1/\ell},
\end{align}
\eqref{First-moment-martingale} and \eqref{Second-moment-martingale}, we get
\begin{align}
	&\sum_{n=1}^\infty \E_x\left( \left|M_{n+1}^{(2k+1,\theta)}- M_{n}^{(2k+1,\theta)}\right|\right)\nonumber\\
	& \leq 2\sum_{n=1}^\infty \E_x\left( \left|M_{n+1}^{(2k+1,\theta)}- M_{n+1}^{(2k+1,\theta, B)}\right|\right) +
  \sum_{n=1}^\infty \mathbb{E}_x\left(
  \left|M_{n+1}^{(2k+1,\theta, B)}- \mathbb{E}_x\left(M_{n+1}^{(2k+1,\theta, B)}\big|\mathcal{F}_n\right)\right|^\ell  \right)^{1/\ell}\nonumber\\
	&\lesssim \sum_{n=1}^\infty \frac{n^k}{n^{1+\lambda}}+ \sum_{n=1}^\infty
	n^{2k +1}e^{-c_1  n},
\end{align}
which is finite since $\lambda >k$.
Therefore,
 $M_n^{(2k+1,\theta)}$
 converges to a limit $M_\infty^{(2k+1,\theta)}$ $\mathbb{P}_x$-a.s. and in $L^1(\P_x)$.

{\bf Step 3:} In this step, we prove the second assertion of (ii).
For any $\eta\in (0, \lambda - k)$,  by \eqref{First-moment-martingale} and \eqref{Second-moment-martingale},
\begin{align}
	& \sum_{n=1}^\infty n^{\lambda-k-\eta}\E_x\left( \left|M_{n+1}^{(2k+1,\theta)}- M_{n}^{(2k+1,\theta)}\right|\right)\lesssim  \sum_{n=1}^\infty n^{\lambda-k-\eta} \frac{n^k}{n^{1+\lambda}}+ \sum_{n=1}^\infty n^{\lambda-k-\eta}
		 n^{2k +1}e^{-c_1  n}<\infty.
\end{align}
Thus,
$n^{\lambda -k-\eta}\left(M_n^{(2k+1,\theta)}-M_\infty^{(2k+1,\theta)}\right)\stackrel{n\to\infty}{\longrightarrow}0,\mathbb{P}_x$-a.s. (see for example \cite[Lemma 2]{Asmussen76}).
For $s\in[n, n+1]$, by Doob's inequality, for any $\varepsilon>0$,
\begin{align}
	&\sum_{n=1}^\infty \mathbb{P}_x\left( n^{\lambda -k-\eta} \sup_{n\leq s\leq n+1} \left| M_s^{(2k+1, \theta)} -M_n^{(2k+1,\theta)}\right| > \varepsilon \right)\nonumber\\
	&\leq \frac{1}{\varepsilon} \sum_{n=1}^\infty n^{\lambda-k -\eta } \mathbb{E}_x\left(\left|M_{n+1}^{(2k+1,\theta)}-M_{n}^{(2k+1,\theta)}\right|\right) <\infty,
\end{align}
which implies that $n^{\lambda-k-\eta} \sup_{n\leq s\leq n+1} \left| M_s^{(2k+1,\theta)} -M_n^{(2k+1,\theta)}\right|  \stackrel{n\to\infty}{\longrightarrow} 0, \mathbb{P}_x$-a.s.
Therefore, we have $\mathbb{P}_x$-almost surely,
\begin{align}
	&\sup_{n\leq s\leq n+1} s^{ \lambda-k -\eta } \left|M_s^{(2k+1,\theta)} - M_\infty^{(2k+1,\theta)} \right| \leq (n+1)^{\lambda -k -\eta} \sup_{n\leq s\leq n+1} \left| M_s^{(2k+1,\theta)} -M_n^{(2k+1,\theta)}\right| \\
	&\quad\quad\quad + (n+1)^{\lambda -k-\eta} \left|M_n^{(2k+1,\theta)}-M_\infty^{(2k+1, \theta)}\right|\stackrel{n\to\infty}{\longrightarrow}0,
\end{align}
which completes the proof of (ii).
\hfill$\Box$

\subsection{Asymptotic expansions along discrete time}

\begin{lemma}\label{lemma2}
Assume $x>0$ and  $\theta \in [0,\sqrt{2})$.
Let $A_n\subset (0,\infty)$ be a family of Borel sets such that
$$
\mbox {either } \quad
 \sup_n \sup\{y: y\in A_n\}<\infty \quad \mbox{ or }\quad  \inf_n \sup\{y: y\in A_n\}=\infty.
$$
For any given $m\in \N$ and $\kappa>1$, if \eqref{LlogL} holds for some
$\lambda>2m+2\kappa+1$, then
	\begin{align}
		 r_n^m 	\frac{Z_{r_n}^{(0,\infty)  } (A_n)-\mathbb{E}_x\left(Z_{r_n}^{(0,\infty)} (A_n)\Big| \mathcal{F}_{\sqrt{r_n}}\right)}{r_n^{-b_\theta}e^{(1-\frac{\theta^2}{2})r_n}}=0,\quad
		 \P_x\mbox{-a.s.}
	\end{align}
where
$$
b_\theta:=\begin{cases}3/2, &
	 \theta=0 \mbox{ and }
	 \sup_n \sup\{y: y\in A_n\}<\infty
	 \mbox{ or } \theta \in (0,\sqrt{2}) ;\\
1/2, & \theta=0 \mbox{ and }
\inf_n \sup\{y: y\in A_n\}=\infty.
\end{cases}
$$
    In particular, for any Borel set $A\subset (0,\infty)$,  as $n\to\infty$,
    	\begin{align}
    	r_n^m 	\frac{Z_{r_n}^{(0,\infty)  } (A)-\mathbb{E}_x\left(Z_{r_n}^{(0,\infty)} (A)\Big| \mathcal{F}_{\sqrt{r_n}}\right)}{r_n^{-b_\theta }e^{(1-\frac{\theta^2}{2})r_n}}=0,\quad
	\P_x\mbox{-a.s.}
    \end{align}
\end{lemma}
\textbf{Proof: }
Suppose $m\in \N$, $\kappa>1$ and that \eqref{LlogL} holds for some
$\lambda>2m+2\kappa+1$.
We divide the proof into three steps.
In Step 1, we define a truncated process $Z_{r_n}^{(0,\infty), G} (A_n)$ and give a first-moment estimate for $Z_{r_n}^{(0,\infty)} (A_n)-Z_{r_n}^{(0,\infty), G} (A_n)$, see \eqref{step_11} below. In Step 2, we
bound the $\ell$-th moment of $Z_{r_n}^{(0,\infty), G} (A_n) - \E_x\left(Z_{r_n}^{(0,\infty), G} (A_n) \big| \mathcal{F}_{\sqrt{r_n}}\right)$ for appropriate $\ell\in(1,2)$,
see \eqref{step_10} below; In Step 3, we combine the results obtained in Step 1 and Step 2
to get the assertion of the proposition.

{\bf Step 1:}
Recall that $v< u$ and $v\leq u$ mean that  $v$ is an ancestor of $u$  and that $v=u$ or $v< u$ respectively.  For $u\in N(r_n)$, define $G_{n,u}$ to be the event that, for all $v<u$ with death time $d_v\in (\sqrt{r_n}, r_n)$, it holds that
$O_v \leq e^{c_0\sqrt{r_n}}$,
where $c_0>0$ is a small constant to be determined later. Define
\begin{align}
	Z_{r_n}^{(0,\infty), G} (A_n):= \sum_{u\in N(r_n)} 1_{\{\min_{s\leq r_n} X_u(s)> 0\}} 	1_{A_n}(X_u(r_n))1_{G_{n,u}}.
\end{align}
 By the branching property, it holds that
\begin{align}\label{Branching-property-1}
	&Z_{r_n}^{(0,\infty), G} (A_n)\nonumber\\
	&= \sum_{v\in N(\sqrt{r_n})} 1_{\left\{\min_{s\leq \sqrt{r_n}} X_v(s)>0 \right\}} \sum_{u\in N(r_n): v\leq u} 1_{\left\{\min_{\sqrt{r_n}<s\leq r_n} X_u(s)>0\right\}}1_{A_n}(X_u(r_n))1_{G_{n,u}}.
\end{align}
Therefore, by the Markov property,
\begin{align}\label{step_0}
	& \frac{r_n^{m+b_\theta}}{ e^{(1-\frac{\theta^2}{2})r_n}} \mathbb{E}_x \left(
	Z_{r_n}^{(0,\infty)} (A_n)- Z_{r_n}^{(0,\infty), G} (A_n)
	\big| \mathcal{F}_{\sqrt{r_n}}\right)\nonumber \\
	& = \frac{r_n^{m+b_\theta}}{ e^{(1-\frac{\theta^2}{2})r_n}}  \sum_{v\in N(\sqrt{r_n})} 1_{\left\{\min_{s\leq \sqrt{r_n}} X_v(s)>0 \right\}} \nonumber\\
	&\quad\quad \times  \mathbb{E}_{X_v(\sqrt{r_n})}\left( \sum_{u\in N(r_n-\sqrt{r_n})}1_{\left\{\min_{s\leq r_n -\sqrt{r_n}} X_u(s)>0 \right\}} 1_{A_n}\left(X_u(r_n-\sqrt{r_n})\right)1_{D_{n,u}^c} \right)
	\nonumber\\ &=:  \frac{r_n^{m+b_\theta}}{ e^{(1-\frac{\theta^2}{2})r_n}}  \sum_{v\in N(\sqrt{r_n})} 1_{\left\{\min_{s\leq \sqrt{r_n}} X_v(s)>0 \right\}}
	 T_{A_n}(X_v(\sqrt{r_n}), r_n-\sqrt{r_n}),
\end{align}
where for $u\in N(r_n-\sqrt{r_n})$, $D_{n,u}$ denotes the event that,
for all $w<u$ with $d_w<r_n-\sqrt{r_n}$,
 it holds that $O_w\leq e^{c_0 \sqrt{r_n}}$.
 Let $d_i$ be the $i$-th splitting time of the spine and $O_i$ be the number of children produced by the spine at time $d_i$.
Define $D_{n, \xi_{r_n-\sqrt{r_n}}}$ to be the event that, for all $i$ with $d_i<r_n -\sqrt{r_n}$, it holds that
$O_i\leq e^{c_0 \sqrt{r_n}}$.
 By Lemma \ref{General-many-to-one}, we get that
\begin{align}\label{step_1}
	&T_{A_n}(X_v(\sqrt{r_n}), r_n-\sqrt{r_n})\nonumber\\
	& = e^{r_n-\sqrt{r_n}} 	\mathbf{E}_{X_v(\sqrt{r_n})}
	\left(1_{\left\{\min_{s\leq r_n -\sqrt{r_n}} X_{\xi}(s)>0  \right\}}1_{A_n}(X_\xi(r_n-\sqrt{r_n}))1_{D_{n,\xi_{r_n-\sqrt{r_n}}}^c} \right).
\end{align}
Note that given $X_\xi$, $\left\{d_i: i\geq 1 \right\}$ are the atoms for a Poisson point process with rate
2, $\{O_i: i\geq 1\}$ are iid with common law $\widehat{L}$ given by
$\mathbf{P}_x(\widehat{L}=k)= kp_k/2$,
and that
$\left\{d_i: i\geq 1 \right\}$ and $\{O_i: i\geq 1\}$ are independent. By \eqref{LlogL}, we conclude that
\begin{align}\label{step_2}
	&\mathbf{E}_{X_v(\sqrt{r_n})}
	\left(1_{D_{n,\xi_{r_n-\sqrt{r_n}}}^c} \big| X_\xi(s): s\geq 0\right)\leq
	\mathbf{E}_{X_v(\sqrt{r_n})}
	\left(\sum_{i: d_i <r_n-\sqrt{r_n}}
	1_{\left\{O_i> e^{c_0 \sqrt{r_n}}\right\}}
	\big| X_\xi(s): s\geq 0\right)\nonumber\\
	& =2 \int_0^{r_n-\sqrt{r_n}} \mathbf{P}_{X_v(\sqrt{r_n})}
          \left(\widehat{L}> e^{c_0 \sqrt{r_n}}\right)\mathrm{d}s \leq 2r_n \frac{
          \mathbf{E}_x
          \left(\log_+^{1+\lambda}L\right)}{(c_0\sqrt{r_n})^{1+\lambda}}\lesssim \frac{1}{r_n^{(\lambda-1)/2}}.
\end{align}
Plugging \eqref{step_2} into \eqref{step_1},  we get
\begin{align}\label{step_3}
	T_{A_n}(X_v(\sqrt{r_n}), r_n-\sqrt{r_n}).
	& \lesssim
	  \frac{e^{r_n-\sqrt{r_n}}}{r_n^{(\lambda-1)/2}}	\mathbf{P}_{X_v(\sqrt{r_n})}
	\left( \min_{s\leq r_n -\sqrt{r_n}} X_\xi(s)>0, X_{r_n-\sqrt{r_n}} \in A_n \right).
\end{align}
Since, under $\mathbf{P}_{X_v(\sqrt{r_n})}$,
$X_\xi(t)$ is a standard Brownian motion with drift $-\theta$. By Lemma \ref{Sharp-Probability} (i) and (ii) with $B=A_n$, we get that
\begin{align}\label{step_4}
	&\mathbf{P}_{X_v(\sqrt{r_n})}
	\left( \min_{s\leq r_n -\sqrt{r_n}} X_\xi(s)>0, X_{r_n-\sqrt{r_n}} \in A_n  \right)\nonumber\\
	&\lesssim \frac{X_v(\sqrt{r_n})e^{\theta X_v(\sqrt{r_n})}}{(r_n-\sqrt{r_n})^{b_\theta} e^{\frac{\theta^2}{2}(r_n-\sqrt{r_n})}} \lesssim  \frac{X_v(\sqrt{r_n})e^{\theta X_v(\sqrt{r_n})}}{r_n^{b_\theta} e^{\frac{\theta^2}{2}(r_n-\sqrt{r_n})}}.
\end{align}
Combining \eqref{step_0}, \eqref{step_3} and \eqref{step_4}, we get that
\begin{align}
	& \frac{r_n^{m+b_\theta}}{ e^{(1-\frac{\theta^2}{2})r_n}} \mathbb{E}_x \left(
	 Z_{r_n}^{(0,\infty)} (A_n)- Z_{r_n}^{(0,\infty), G} (A_n)
	\big| \mathcal{F}_{\sqrt{r_n}}\right)\nonumber \\
	& \lesssim  \frac{r_n^{m+b_\theta}}{ e^{(1-\frac{\theta^2}{2})r_n}}  \sum_{v\in N(\sqrt{r_n})} 1_{\left\{\min_{s\leq \sqrt{r_n}} X_v(s)>0 \right\}}
  	  \frac{e^{r_n-\sqrt{r_n}}}{r_n^{(\lambda-1)/2}}
	\frac{X_v(\sqrt{r_n})e^{\theta X_v(\sqrt{r_n})}}{r_n^{b_\theta} e^{\frac{\theta^2}{2}(r_n-\sqrt{r_n})}}.\\
	 & = \frac{1}{r_n^{(\lambda-1-2m)/2}}
	e^{-(1-\frac{\theta^2}{2})\sqrt{r_n}}\sum_{v\in N(\sqrt{r_n})}1_{\left\{\min_{s\leq \sqrt{r_n}} X_v(s)>0 \right\}} X_v(\sqrt{r_n}) e^{\theta X_v(\sqrt{r_n})} =
	 \frac{1}{r_n^{(\lambda-1-2m)/2}}
	M_{\sqrt{r_n}}^{(1,\theta)},
\end{align}
with $M_{\sqrt{r_n}}^{(1,\theta)}$ given in \eqref{Martingale}.
Now taking expectation with respect to $\P_x$,
we get that
\begin{align}\label{step_11}
		& \frac{r_n^{m+b_\theta}}{ e^{(1-\frac{\theta^2}{2})r_n}} \mathbb{E}_x \left(
		Z_{r_n}^{(0,\infty)} (A_n)- Z_{r_n}^{(0,\infty), G} (A_n)
		\right)
		 \lesssim \frac{1}{r_n^{(\lambda-1-2m)/2}} .
\end{align}

{\bf Step 2:}
By \eqref{Branching-property-1}
 and the branching property,
\begin{align}
	&Z_{r_n}^{(0,\infty), G} (A_n)- \mathbb{E}_x \left(Z_{r_n}^{(0,\infty), G} (A_n)\big| \mathcal{F}_{\sqrt{r_n}} \right)=: \sum_{v\in N(\sqrt{r_n})} 1_{\left\{\min_{s\leq \sqrt{r_n}} X_v(s)>0 \right\}} J_{n,u},
\end{align}
where conditioned on $\mathcal{F}_{\sqrt{r_n}}$,  $\left\{J_{n,v}: v\in N(\sqrt{r_n})\right\}$ are centered independent random variables defined by
\begin{align}
	J_{n,v} =&\sum_{u\in N(r_n): v\leq u} 1_{\left\{\min_{\sqrt{r_n}<s\leq r_n} X_u(s)>0\right\}}1_{A_n}(X_u(r_n))1_{G_{n,u}} \nonumber\\
	&  -\mathbb{E}_x\left(\sum_{u\in N(r_n): v\leq u} 1_{\left\{\min_{\sqrt{r_n}<s\leq r_n} X_u(s)>0\right\}}1_{A_n}(X_u(r_n))1_{G_{n,u}}\Big| \mathcal{F}_{\sqrt{r_n}}\right) .
\end{align}
 Thus, by Lemma \ref{Useful-Ineq-1},  for any fixed $\ell$
 with $1< \ell <\min\{2, 2/\theta^2\}$ and  $(\ell-1)\frac{\theta^2}{2}< 1-\frac{\theta^2}{2}$,
\begin{align}\label{step_8}
&\mathbb{E}_x \left(\left|Z_{r_n}^{(0,\infty), G} (A_n)- \mathbb{E}_x \left(Z_{r_n}^{(0,\infty), G} (A_n)\big| \mathcal{F}_{\sqrt{r_n}} \right)\right|^\ell \bigg| \mathcal{F}_{\sqrt{r_n}}\right) \nonumber \\
&\leq 2\sum_{v\in N(\sqrt{r_n})} 1_{\left\{\min_{s\leq \sqrt{r_n}} X_v(s)>0 \right\}}  \mathbb{E}_x\left( \left|J_{n,v}\right|^\ell  \big| \mathcal{F}_{\sqrt{r_n}} \right) \lesssim \sum_{v\in N(\sqrt{r_n})} M_{n,v},
\end{align}
where for each $v\in N(\sqrt{r_n})$,
\begin{align}
	M_{n,v}:= \mathbb{E}_{X_v(\sqrt{r_n})}\left(\left(\sum_{u\in N(r_n-\sqrt{r_n})}1_{\left\{\min_{s\leq r_n- \sqrt{r_n}} X_u(s)>0 \right\}} 1_{A_n}\left(X_u(r_n-\sqrt{r_n})\right)1_{D_{n,u}} \right)^\ell \right).
\end{align}
Set $V_n:= \sum_{u\in N(r_n-\sqrt{r_n})}1_{\left\{\min_{s\leq r_n- \sqrt{r_n}} X_u(s)>0 \right\}} 1_{A_n}\left(X_u(r_n-\sqrt{r_n})\right)1_{D_{n,u}} $. By Lemma \ref{General-many-to-one} and the fact that $A_n \subset (0,\infty), V_n \leq \sum_{u\in N(r_n-\sqrt{r_n})} 1_{(0,\infty)}(X_u(r_n-\sqrt{r_n}))$, we have
\begin{align}\label{step_7}
	& M_{n,v}= e^{r_n-\sqrt{r_n}}	\mathbf{E}_{X_v(\sqrt{r_n})}
	\left(1_{\left\{\min_{s\leq r_n- \sqrt{r_n}} X_\xi(s)>0 \right\}} 1_{A_n}\left(X_\xi(r_n-\sqrt{r_n})\right)1_{D_{n,\xi_{r_n-\sqrt{r_n}} }} V_n^{\ell-1} \right)\nonumber\\
& \leq e^{r_n-\sqrt{r_n}} \mathbf{E}_{X_v(\sqrt{r_n})}
\left(1_{D_{n,\xi_{r_n-\sqrt{r_n}} }} 1_{\{X_\xi(r_n-\sqrt{r_n})> 0\}}\left(\sum_{u\in N(r_n-\sqrt{r_n})} 1_{\left\{ X_u(r_n-\sqrt{r_n})>0 \right\}}  \right)^{\ell-1}\right).
\end{align}
Given $X_\xi, d_i $ and $O_i$,
by the Markov property and the inequality $(\sum_{i=1}^n x_i )^p \leq \sum_{i=1}^n x_i^p$ for all $x_i\geq 0$ and $p\in (0,1)$, it holds that
\begin{align}\label{step_5}
&1_{\left\{D_{n,\xi_{r_n-\sqrt{r_n}} }\right\}}
\mathbf{E}_{X_v(\sqrt{r_n})}
\left( \left(\sum_{u\in N(r_n-\sqrt{r_n})} 1_{\left\{ X_u(r_n-\sqrt{r_n})>0 \right\}}  \right)^{\ell-1}\big| X_\xi, d_i, O_i: i\geq 1\right)\nonumber \\
& \leq 1_{\left\{D_{n,\xi_{r_n-\sqrt{r_n}} }\right\}} \sum_{i: d_i < r_n-\sqrt{r_n}} (O_i -1)\mathbb{E}_{X_\xi(d_i)} \left( \left(\sum_{u\in N(z)} 1_{\left\{ X_u(r_n-\sqrt{r_n})>0 \right\}}\right)^{\ell-1} \right)\big|_{z=r_n-\sqrt{r_n}-d_i} \nonumber \\
& \leq
e^{c_0 \sqrt{r_n}} \sum_{i: d_i < r_n-\sqrt{r_n}}  \mathbb{E}_{X_\xi(d_i)} \left( \sum_{u\in N(z)} 1_{\left\{ X_u(r_n-\sqrt{r_n})>0 \right\}}  \right)^{\ell-1}\big|_{z=r_n-\sqrt{r_n}-d_i},
\end{align}
where in the last inequality, we
used the  fact that, on the event $D_{n, \xi_{r_n-\sqrt{r_n}}}$, $O_i-1\leq e^{c_0 \sqrt{r_n}}$ and the fact that ${\rm E}(|X|^p)\leq {\rm E} (|X|)^p$ for $p\in (0,1)$.
Note that by Lemma \ref{General-many-to-one} and Lemma \ref{Sharp-Probability} (iii),
\[
\mathbb{E}_x \left( \sum_{u\in N(z)} 1_{\left\{ X_u(r_n-\sqrt{r_n})>0 \right\}}  \right)\leq e^{(1-\frac{\theta^2}{2})z}e^{\theta x}.
\]
Using the fact that $d_i$ are the atoms of a Poisson process with rate $2$, taking expectation with respect to
$\mathbf{P}_{X_v(\sqrt{r_n})} \left(\cdot |X_\xi\right)$
in \eqref{step_5}, we get that
\begin{align}\label{step_6}
	&\mathbf{E}_{X_v(\sqrt{r_n})}
	\left(1_{\left\{D_{n,\xi_{r_n-\sqrt{r_n}} }\right\}}\left(\sum_{u\in N(r_n-\sqrt{r_n})} 1_{\left\{ X_u(r_n-\sqrt{r_n})>0 \right\}}   \right)^{\ell-1} \big| X_\xi \right)\nonumber\\
	 &\leq
	 2e^{c_0 \sqrt{r_n}}
	  \int_0^{r_n-\sqrt{r_n}}
	  e^{\theta (\ell-1)X_\xi(s)} e^{(1-\frac{\theta^2}{2})(\ell-1)\left(r_n-\sqrt{r_n}-s\right)}\mathrm{d}s.
\end{align}
Combining \eqref{step_7} and \eqref{step_6},
noting that $X_\xi(s)$ under $\mathbf{P}_x$
is a standard Brownian motion with drift $-\theta$,
and applying Lemma \ref{Sharp-Probability}(iii), we conclude that
\begin{align}
	&M_{n,v}\lesssim e^{r_n-\sqrt{r_n}}	\mathbf{E}_{X_v(\sqrt{r_n})}
	\left( e^{c_0\sqrt{ r_n}}		1_{\{X_\xi(r_n-\sqrt{r_n})> 0\}}  \int_0^{r_n-\sqrt{r_n}} e^{\theta (\ell-1)X_\xi(s)} e^{(1-\frac{\theta^2}{2})(\ell-1)\left(r_n-\sqrt{r_n}-s\right)}\mathrm{d}s\right)  \nonumber\\
	& = e^{r_n-\sqrt{r_n}}
	e^{c_0 \sqrt{r_n}}	 \int_0^{r_n-\sqrt{r_n}}	  \Pi_{X_v(\sqrt{r_n})}^{-\theta} \left(e^{\theta (\ell-1)B_s} \Pi_{B_s}^{-\theta}	\left(B_{r_n-\sqrt{r_n} -s}>0 \right)  \right)
	  e^{(1-\frac{\theta^2}{2})(\ell-1)\left(r_n-\sqrt{r_n}-s\right)}\mathrm{d}s \\
 & \leq  e^{r_n-\sqrt{r_n}} e^{c_0 \sqrt{r_n}} \int_0^{r_n-\sqrt{r_n}}
  \Pi_{X_v(\sqrt{r_n})}^{-\theta} \left(e^{\ell \theta  B_s} \right) e^{-\frac{\theta^2}{2}(r_n-\sqrt{r_n}-s)}	  e^{(1-\frac{\theta^2}{2})(\ell-1)\left(r_n-\sqrt{r_n}-s\right)}\mathrm{d}s.
\end{align}
Using elementary calculus in the last integral, we get
\begin{align}
	&M_{n,v}\lesssim e^{\ell \theta X_v(\sqrt{r_n})} e^{\ell (1-\frac{\theta^2}{2})(r_n-\sqrt{r_n})}		e^{c_0 \sqrt{r_n}}			 \int_0^{r_n-\sqrt{r_n}}   e^{-s(\ell-1)(1-\frac{\theta^2}{2}-(\ell-1)\frac{\theta^2}{2})}\mathrm{d}s\\
	& \lesssim e^{\ell \theta X_v(\sqrt{r_n})} e^{\ell (1-\frac{\theta^2}{2})(r_n-\sqrt{r_n})}		e^{c_0 \sqrt{r_n}}.
\end{align}
Plugging this upper-bound into \eqref{step_8},  we obtain that
\begin{align}\label{step_9}
	&\frac{r_n^{\ell (m+b_\theta)}}{ e^{\ell (1-\frac{\theta^2}{2})r_n}}\mathbb{E}_x \left(\left|Z_{r_n}^{(0,\infty), G} (A_n)- \mathbb{E}_x \left(Z_{r_n}^{(0,\infty), G} (A_n)\big| \mathcal{F}_{\sqrt{r_n}} \right)\right|^\ell \bigg| \mathcal{F}_{\sqrt{r_n}}\right) \nonumber\\
    & \lesssim  \frac{r_n^{\ell (m+b_\theta)}}{ e^{\ell (1-\frac{\theta^2}{2})r_n}}\sum_{v\in N(\sqrt{r_n})} e^{\theta \ell X_v(\sqrt{r_n})} e^{\ell (1-\frac{\theta^2}{2})(r_n-\sqrt{r_n})}		e^{c_0 \sqrt{r_n}}
    \nonumber\\ & =\frac{ r_n^{\ell(m+b_\theta)} e^{c_0\sqrt{r_n}}}{e^{\ell (1-\frac{\theta^2}{2})\sqrt{r_n}}}    \sum_{v\in N(\sqrt{r_n})} e^{\theta \ell X_v(\sqrt{r_n})}.
\end{align}
 Taking expectation in \eqref{step_9} with respect to $\P_x$, and using  Lemma \ref{General-many-to-one},
 we get
\begin{align}\label{step_10}
& \frac{r_n^{\ell (m+b_\theta)}}{ e^{\ell (1-\frac{\theta^2}{2})r_n}}\mathbb{E}_x \left(\left|Z_{r_n}^{(0,\infty), G} (A_n)- \mathbb{E}_x \left(Z_{r_n}^{(0,\infty), G} (A_n)\big| \mathcal{F}_{\sqrt{r_n}} \right)\right|^\ell \right) \nonumber\\
    & \lesssim \frac{ r_n^{\ell(m+b_\theta)} e^{c_0\sqrt{r_n}}}{e^{\ell (1-\frac{\theta^2}{2})\sqrt{r_n}}}    \E_x \left(\sum_{v\in N(\sqrt{r_n})} e^{\theta \ell X_v(\sqrt{r_n})} \right)  = \frac{ r_n^{\ell(m+b_\theta)} e^{c_0\sqrt{r_n}}}{e^{\ell (1-\frac{\theta^2}{2})\sqrt{r_n}}}  e^{\sqrt{r_n}}  \mathbf{E}_x \left(e^{\theta \ell X_\xi(\sqrt{r_n})} \right)
\nonumber\\ & = r_n^{\ell (m+b_\theta)}
e^{\theta \ell x}
e^{-((\ell-1)(1-\frac{\theta^2\ell}{2})-c_0 )\sqrt{r_n}}.
\end{align}

{\bf Step 3:}
Fix $c_0\in (0, (\ell-1)(1-\frac{\theta^2\ell}{2}))$
and set $c_1:= ((\ell-1)(1-\frac{\theta^2\ell}{2})-c_0 )/\ell>0$.
Using \eqref{New-step_2}
 with $X= Z_{r_n}^{(0,\infty)}(A_n)$ and $Y=Z_{r_n}^{(0,\infty), G}(A_n),$
  we have
\begin{align}
		&\mathbb{E}_x\left(\frac{r_n^m}{r_n^{-b_\theta}e^{(1-\frac{\theta^2}{2})r_n}} \left|Z_{r_n}^{(0,\infty)} (A_n)- \mathbb{E}_x \left(Z_{r_n}^{(0,\infty)} (A_n)\big| \mathcal{F}_{\sqrt{r_n}} \right)\right|\right)\nonumber\\
		& \leq   2 \frac{r_n^{m+b_\theta}}{ e^{(1-\frac{\theta^2}{2})r_n}} \mathbb{E}_x \left(
		Z_{r_n}^{(0,\infty)} (A_n)- Z_{r_n}^{(0,\infty), G} (A_n)\right) \nonumber\\
		&\quad + \left(\frac{r_n^{\ell (m+b_\theta)}}{ e^{\ell (1-\frac{\theta^2}{2})r_n}}\mathbb{E}_x \left(\left|Z_{r_n}^{(0,\infty), G} (A_n)- \mathbb{E}_x \left(Z_{r_n}^{(0,\infty), G} (A_n)\big| \mathcal{F}_{\sqrt{r_n}} \right)\right|^\ell \right) \right)^{1/\ell} \nonumber\\
		& \lesssim
		\frac{1}{r_n^{(\lambda-1-2m)/2}} + r_n^{(m+b_\theta)}  e^{-c_1 \sqrt{r_n}} = \frac{1}{n^{(\lambda-1-2m)/(2\kappa)}}+ n^{(m+b_\theta)/\kappa}e^{-c_1n^{1/(2\kappa)}},
\end{align}
where in the last inequality, we used  \eqref{step_11} and \eqref{step_10}.
Since $\lambda> 2m+2\kappa+1$,
we conclude that for any $\varepsilon>0$,
\begin{align}
	&\sum_{n=1}^\infty \mathbb{P}_x\left(  \frac{r_n^m}{r_n^{-b_\theta}e^{(1-\frac{\theta^2}{2})r_n}} \left|Z_{r_n}^{(0,\infty)} (A_n)- \mathbb{E}_x \left(Z_{r_n}^{(0,\infty)} (A_n)\big| \mathcal{F}_{\sqrt{r_n}} \right)\right| >\varepsilon \right)\\
	& \lesssim\sum_{n=1}^\infty \left(
	\frac{1}{n^{(\lambda-1-2m)/(2\kappa)}}+ n^{(m+b_\theta)/\kappa}e^{-c_1n^{1/(2\kappa)}}\right)<\infty,
\end{align}
which completes the proof of the Lemma.
\hfill$\Box$

\begin{prop}\label{Discrete-time-expansion}
Let $x>0$.	For any given $m\in \N$ and
$\kappa>2m+1$,
if \eqref{LlogL} holds for some
$\lambda>2m+2\kappa +2$,
then
	(i)
	for any $\theta\in (0,\sqrt{2})$, $x>0$ and Borel set $A\subset(0,\infty)$, $\P_x$-almost surely, as $n\to\infty$,
	\begin{align}
		&\frac{Z_{r_n}^{(0,\infty)} (A)}{r_n^{-3/2}e^{(1-\frac{\theta^2}{2})r_n}}\nonumber\\
		&=- \sqrt{\frac{2}{\pi}}\sum_{\ell=0}^m \frac{H_{2\ell+2}(0)}{r_n^{\ell}} \sum_{k=0}^{\ell} \frac{M_\infty^{(2k+1,\theta)} }{(2k+1)! (2\ell -2k+1)!}
		\int_A z^{2\ell -2k+1} e^{-\theta z} \mathrm{d} z+o(r_n^{-m});
	\end{align}
	(ii)
	 for any bounded Borel set $A\subset (0,\infty)$,
	$\P_x$-almost surely, as $n\to\infty$,
	\begin{align}
		&\frac{Z_{r_n}^{(0,\infty)} (A)}{r_n^{-3/2}e^{r_n}}\nonumber\\
		&=- \sqrt{\frac{2}{\pi}}\sum_{\ell=0}^m \frac{H_{2\ell+2}(0)}{r_n^{\ell}} \sum_{k=0}^{\ell} \frac{M_\infty^{(2k+1,0)} }{(2k+1)! (2\ell -2k+1)!}
		\int_A z^{2\ell -2k+1}  \mathrm{d} z+o(r_n^{-m});
	\end{align}
	(iii)
	for any $a\geq 0$, $\P_x$-almost surely, as $n\to\infty$,
	\begin{align}
		&\frac{Z_{r_n}^{(0,\infty)} ((a,\infty))}{r_n^{-1/2}e^{r_n}}= \sqrt{\frac{2}{\pi}}\sum_{\ell=0}^m \frac{H_{2\ell}(0)}{r_n^{\ell}} \sum_{k=0}^{\ell} \frac{M_\infty^{(2k+1,0)} }{(2k+1)! (2\ell -2k)!} a^{2\ell-2k} +o(r_n^{-m}).
	\end{align}
\end{prop}
\textbf{Proof: }
Suppose $x>0, m\in \N, \kappa>2m+1$ and \eqref{LlogL} holds with $\lambda> 2m+ 2\kappa +2> 2m+\kappa+1$.

 First using Lemma \ref{lemma2} and then the Markov property, we get that
 for any Borel set $A\subset (0,\infty)$,
$\P_x$-almost surely,
\begin{align}\label{step_12}
	& \frac{Z_{r_n}^{(0,\infty)  } (A)}{r_n^{-b_\theta}e^{(1-\frac{\theta^2}{2})r_n}} = o(r_n^{-m})+ \frac{\mathbb{E}_x\left(Z_{r_n}^{(0,\infty)} (A)\Big| \mathcal{F}_{\sqrt{r_n}}\right)}{r_n^{-b_\theta}e^{(1-\frac{\theta^2}{2})r_n}}\nonumber\\
&	= o(r_n^{-m})+ \sum_{v\in N(\sqrt{r_n})} \frac{1_{\{\min_{s\leq \sqrt{r_n}}
		X_v(s)
		>0 \}}}{r_n^{-b_\theta }e^{(1-\frac{\theta^2}{2})r_n}} \mathbb{E}_{X_v(\sqrt{r_n})}\left(Z_{r_n-\sqrt{r_n}}^{(0,\infty)} (A) \right)\nonumber\\
&= o(r_n^{-m})+ \sum_{v\in N(\sqrt{r_n})} \frac{1_{\{\min_{s\leq \sqrt{r_n}} X_v(s)>0\}}}{r_n^{-b_\theta}e^{(1-\frac{\theta^2}{2})r_n}}\nonumber\\
&\quad \quad \times e^{r_n-\sqrt{r_n}}\Pi_{X_v(\sqrt{r_n})}^{-\theta}\left( 1_{\{\min_{s\leq r_n -\sqrt{r_n}} B_s>0\}} 1_{A}(B_{r_n-\sqrt{r_n}}) \right),
\end{align}
where in the last equality we used Lemma \ref{General-many-to-one}.
Let  $K:= 2m/\kappa +3$ and fix a sufficient small $\varepsilon>0$  such that
\begin{align}\label{Small-Varepsilon}
	K(1-\varepsilon)> \frac{2m+1}{\kappa}+2.
\end{align}

\textbf{Step 1:}  In this step, we prove that  $\P_x$-almost surely,
\begin{align}\label{step_40}
	&\frac{Z_{r_n}^{(0,\infty)  } (A)}{r_n^{-b_\theta}e^{(1-\frac{\theta^2}{2})r_n}} = o(r_n^{-m})+ \sum_{v\in N(\sqrt{r_n})} \frac{1_{\{\min_{s\leq \sqrt{r_n}} X_v(s) > 0 \}}}{r_n^{-b_\theta}e^{(1-\frac{\theta^2}{2})r_n}} 1_{\{ X_v(\sqrt{r_n})\leq  \sqrt{K\sqrt{r_n}\log n}\}}\nonumber\\
	&\quad \quad \times e^{r_n-\sqrt{r_n}}\Pi_{X_v(\sqrt{r_n})}^{-\theta}\left( 1_{\{\min_{s\leq r_n -\sqrt{r_n}} B_s>0\}} 1_{A}(B_{r_n-\sqrt{r_n}}) \right).
\end{align}
Using  Lemma \ref{Sharp-Probability} (i) (ii) first, and then Lemma \ref{General-many-to-one}, we get
\begin{align}\label{step_14}
	&\sum_{n=2}^\infty r_n^m \mathbb{E}_{x}\bigg(  \sum_{v\in N(\sqrt{r_n})} \frac{1_{\{\min_{s\leq \sqrt{r_n}} X_v(s)>0\}}}{r_n^{-b_\theta}e^{(1-\frac{\theta^2}{2})r_n}} 1_{\{ X_v(\sqrt{r_n})> \sqrt{K\sqrt{r_n}\log n}\}}\nonumber\\
	&\quad \quad \times e^{r_n-\sqrt{r_n}}\Pi_{X_v(\sqrt{r_n})}^{-\theta}\left( 1_{\{\min_{s\leq r_n -\sqrt{r_n}} B_s>0\}} 1_{A}(B_{r_n-\sqrt{r_n}}) \right)\bigg)\nonumber\\
	& \lesssim\sum_{n=2}^\infty r_n^m \mathbb{E}_{x}\bigg(  \sum_{v\in N(\sqrt{r_n})} \frac{1_{\{\min_{s\leq \sqrt{r_n}} X_v(s)>0\}}}{r_n^{-b_\theta }e^{(1-\frac{\theta^2}{2})r_n}} 1_{\{ X_v(\sqrt{r_n})> \sqrt{K\sqrt{r_n}\log n}\}} \nonumber\\
	&\quad \quad \times X_v(\sqrt{r_n})e^{\theta X_v(\sqrt{r_n})}
	\frac{e^{(1-\frac{\theta^2}{2})(r_n-\sqrt{r_n})}}{(r_n-\sqrt{r_n})^{b_\theta}} \bigg)
	\nonumber\\
	&= \sum_{n=2}^\infty \frac{ r_n^m  r_n^{b_\theta} e^{\frac{\theta^2}{2}\sqrt{r_n}}}{\left(r_n-\sqrt{r_n}\right)^{b_\theta}}
	\mathbf{E}_{x}
	\bigg(  1_{\{\min_{s\leq \sqrt{r_n}} X_\xi (s)>0\}}  1_{\{ X_\xi(\sqrt{r_n})> \sqrt{K\sqrt{r_n}\log n}\}}  X_\xi(\sqrt{r_n})e^{\theta X_\xi(\sqrt{r_n})}\bigg).
\end{align}
Recall that $(X_\xi(t), \mathbf{P}_x)$
is equal in law to a standard Brownian motion with drift $-\theta$.
By Lemma \ref{lemma3} (i),
recalling the choice of $\varepsilon$ in \eqref{Small-Varepsilon},
the left-hand side of \eqref{step_14} is bounded from above by
\begin{align}\label{step_16}
&
	\sum_{n=2}^\infty \frac{ r_n^m  r_n^{b_\theta}}{\left(r_n-\sqrt{r_n}\right)^{b_\theta}} xe^{\theta x}\Pi_x^\uparrow \left(B_{\sqrt{r_n}}> \sqrt{K \sqrt{r_n}\log n} \right)\nonumber\\
	&\leq xe^{\theta x}\sum_{n=2}^\infty \frac{ r_n^m  r_n^{b_\theta}}{\left(r_n-\sqrt{r_n}\right)^{b_\theta}} \frac{1}{n^{K(1-\varepsilon)/2}}\Pi_x^\uparrow \left(e^{(1-\varepsilon)\left(B_{\sqrt{r_n}}\right)^2/(2\sqrt{r_n})} \right),
\end{align}
where in the last inequality we used the Markov inequality.
Using Lemma \ref{lemma3} (ii), we get
\begin{align}\label{step_15}
 & \Pi_x^\uparrow	\left(e^{(1-\varepsilon)\left(B_{\sqrt{r_n}}\right)^2/(2\sqrt{r_n})} \right)= \int_0^\infty e^{(1-\varepsilon)y^2/(2\sqrt{r_n})}
     p_{\sqrt{r_n}}^\uparrow (x, y)\mathrm{d}y\nonumber\\
 & \lesssim \int_0^\infty e^{(1-\varepsilon)y^2/(2\sqrt{r_n})} \frac{y^2}{r_n^{3/4}} e^{-(x-y)^2/(2\sqrt{r_n})}\mathrm{d} y = \int_0^\infty e^{(1-\varepsilon)y^2/2} y^2 e^{-(x r_n^{-1/4}-y)^2/2}\mathrm{d}y \lesssim 1.
\end{align}
Combining \eqref{step_16}, \eqref{step_15} and the fact that $r_n (r_n-\sqrt{r_n})^{-1}\lesssim 1$, we obtain
\begin{align}
	&\sum_{n=2}^\infty r_n^m \mathbb{E}_{x}\bigg(  \sum_{v\in N(\sqrt{r_n})} \frac{1_{\{\min_{s\leq \sqrt{r_n}} X_v(s)>0\}}}{r_n^{-b_\theta}e^{(1-\frac{\theta^2}{2})r_n}} 1_{\{ X_v(\sqrt{r_n})> \sqrt{K\sqrt{r_n}\log n}\}}\nonumber\\
		&\quad \quad \times e^{r_n-\sqrt{r_n}}\Pi_{X_v(\sqrt{r_n})}^{-\theta}\left( 1_{\{\min_{s\leq r_n -\sqrt{r_n}} B_s>0\}} 1_{A}(B_{r_n-\sqrt{r_n}}) \right)\bigg)\nonumber\\
	& \lesssim \sum_{n=2}^\infty  \frac{r_n^m}{n^{K(1-\varepsilon)/2}}= \sum_{n=2}^\infty \frac{n^{m/\kappa}}{n^{K(1-\varepsilon)/2}}<\infty,
\end{align}
which implies that $\P_x$-almost surely,
\begin{align}
& r_n^m  \sum_{v\in N(\sqrt{r_n})} \frac{1_{\{\min_{s\leq \sqrt{r_n}} X_v(s)>0\}}}{r_n^{-b_\theta}e^{(1-\frac{\theta^2}{2})r_n}} 1_{\{ X_v(\sqrt{r_n})> \sqrt{K\sqrt{r_n}\log n}\}}\nonumber\\
&\quad \quad \times e^{r_n-\sqrt{r_n}}\Pi_{X_v(\sqrt{r_n})}^{-\theta}\left( 1_{\{\min_{s\leq r_n -\sqrt{r_n}} B_s>0\}} 1_{A}(B_{r_n-\sqrt{r_n}}) \right) \stackrel{n\to\infty}{\longrightarrow}0,
\end{align}
and thus \eqref{step_40} holds.

\textbf{Step 2:}
In this step, we prove (i) and (ii). Recall that when $\theta \in (0,\sqrt{2})$, or $\theta=0$ and
$\sup\{y: y\in A\}<\infty$, $b_\theta=3/2$.
 In this case,  let
 $J:= 2m+2\kappa +2> 2m+ \frac{K\kappa+1}{2}$.
 By using  Lemma \ref{lemma3} (i)  in the first two equalities below, Lemma \ref{lemma1} in the third, we get
\begin{align}
	&\frac{Z_{r_n}^{(0,\infty)  } (A)}{r_n^{-3/2}e^{(1-\frac{\theta^2}{2})r_n}} \nonumber\\
	    & = o(r_n^{-m})+ \sum_{v\in N(\sqrt{r_n})} \frac{1_{\{\min_{s\leq \sqrt{r_n}} X_v(s)>0 \}}}{r_n^{-3/2}e^{(1-\frac{\theta^2}{2})\sqrt{r_n}}} X_v(\sqrt{r_n})e^{\theta X_v(\sqrt{r_n})} 1_{\{ X_v(\sqrt{r_n})\leq  \sqrt{K\sqrt{r_n}\log n}\}}\nonumber\\
	&\quad\quad \times   \Pi_{X_v(\sqrt{r_n})}^\uparrow\left(\frac{1_{A}(B_{r_n-\sqrt{r_n}})}{B_{r_n-\sqrt{r_n}} e^{\theta B_{r_n-\sqrt{r_n}}}}\right)	\nonumber\\
	& = o(r_n^{-m})+ \sum_{v\in N(\sqrt{r_n})} \frac{1_{\{\min_{s\leq \sqrt{r_n}} X_v(s)>0 \}}}{e^{(1-\frac{\theta^2}{2})\sqrt{r_n}}} e^{\theta X_v(\sqrt{r_n})} 1_{\{ X_v(\sqrt{r_n})\leq  \sqrt{K\sqrt{r_n}\log n}\}}\nonumber\\
	&\quad\quad \times \int_A \frac{r_n^{3/2}}{e^{\theta z} \sqrt{r_n-\sqrt{r_n}}}\left(\phi\left( \frac{z-X_v(\sqrt{r_n})}{\sqrt{r_n-\sqrt{r_n}}}\right) -\phi\left( \frac{z+X_v(\sqrt{r_n})}{\sqrt{r_n-\sqrt{r_n}}}\right)  \right)\mathrm{d}z\nonumber\\
	& =  o(r_n^{-m})+ \sum_{v\in N(\sqrt{r_n})} \frac{1_{\{\min_{s\leq \sqrt{r_n}} X_v(s)>0\}}}{e^{(1-\frac{\theta^2}{2})\sqrt{r_n}}} e^{\theta X_v(\sqrt{r_n})} 1_{\{ X_v(\sqrt{r_n})\leq  \sqrt{K\sqrt{r_n}\log n}\}} \int_A r_n e^{-\theta z} \nonumber\\
	& \times
	\left(2\phi\left(\frac{z}{\sqrt{r_n}}\right)  \sum_{k=0}^J \frac{1}{(2k+1)!} \frac{1}{r_n^{(2k+1)/2}} H_{2k+1}\left(\frac{z}{\sqrt{r_n}}\right)  r_n^{(2k+1)/4}H_{2k+1}\left(\frac{X_v(\sqrt{r_n})}{r_n^{1/4}}\right) + \varepsilon_{m,v, z, n}\right)\mathrm{d} z,
\end{align}
where the error term $\varepsilon_{m,v, z, n}$ satisfies that
\[
r_n^{m+1} \sup_{z>0}\sup_{v\in N(\sqrt{r_n})} \left| \varepsilon_{m,v, z, n}\right|1_{\{ X_v(\sqrt{r_n})\leq  \sqrt{K\sqrt{r_n}\log n}\}}\stackrel{n\to\infty}{\longrightarrow}0,\quad
\P_x\mbox{-a.s.}
\]
Noticing that $e^{-(1-\frac{\theta^2}{2})t}\sum_{v\in N(t)} e^{\theta X_v(t)}$ is a non-negative martingale, and that $\int_Ae^{-\theta z}\mathrm{d}z<\infty $
when $\theta>0$, or  $\theta=0$ and
$\sup\{y: y\in A\}<\infty$, we get
\begin{align}
	& r_n^m  \sum_{v\in N(\sqrt{r_n})} \frac{1_{\{\min_{s\leq \sqrt{r_n}} X_v(s)>0\}}}{e^{(1-\frac{\theta^2}{2})\sqrt{r_n}}} e^{\theta X_v(\sqrt{r_n})} 1_{\{ X_v(\sqrt{r_n})\leq  \sqrt{K\sqrt{r_n}\log n}\}} \int_A r_n e^{-\theta z} \left| \varepsilon_{m,v, z, n}\right| \mathrm{d} z\nonumber\\
	& \leq  \Big(e^{-(1-\frac{\theta^2}{2})\sqrt{r_n}}\sum_{v\in N(\sqrt{r_n})} e^{\theta X_v(\sqrt{r_n})}\Big) \left(\int_A e^{-\theta z}\mathrm{d}z\right)\nonumber\\
	& \quad\quad\times r_n^{m+1} \sup_{z>0}\sup_{v\in N(\sqrt{r_n})} \left| \varepsilon_{m,v, z, n}\right|1_{\{ X_v(\sqrt{r_n})\leq  \sqrt{K\sqrt{r_n}\log n}\}}\stackrel{n\to\infty}{\longrightarrow}0,\quad \mbox{a.s.}
\end{align}
Therefore,  $\P_x$-almost surely,
\begin{align}\label{step_19}
	&\frac{Z_{r_n}^{(0,\infty)  } (A)}{r_n^{-3/2}e^{(1-\frac{\theta^2}{2})r_n}} \nonumber\\
	&=o(r_n^{-m})+ \sum_{v\in N(\sqrt{r_n})} \frac{1_{\{\min_{s\leq \sqrt{r_n}} X_v(s)>0\}}}{e^{(1-\frac{\theta^2}{2})\sqrt{r_n}}} e^{\theta X_v(\sqrt{r_n})} 1_{\{ X_v(\sqrt{r_n})\leq  \sqrt{K\sqrt{r_n}\log n}\}} \int_A r_n e^{-\theta z} \nonumber\\
	& \times
	2\phi\left(\frac{z}{\sqrt{r_n}}\right)  \sum_{k=0}^J \frac{1}{(2k+1)!} \frac{1}{r_n^{(2k+1)/2}} H_{2k+1}\left(\frac{z}{\sqrt{r_n}}\right)  r_n^{(2k+1)/4}H_{2k+1}\left(\frac{X_v(\sqrt{r_n})}{r_n^{1/4}}\right) \mathrm{d} z\nonumber\\
	& =o(r_n^{-m})+ 2\sum_{k=0}^J  \sum_{v\in N(\sqrt{r_n})} \frac{1_{\{\min_{s\leq \sqrt{r_n}} X_v(s)>0\}}}{e^{(1-\frac{\theta^2}{2})\sqrt{r_n}}} e^{\theta X_v(\sqrt{r_n})}   r_n^{(2k+1)/4}H_{2k+1}\left(\frac{X_v(\sqrt{r_n})}{r_n^{1/4}}\right)  \nonumber\\
	& \times 1_{\{ X_v(\sqrt{r_n})\leq  \sqrt{K\sqrt{r_n}\log n}\}}
	\int_A e^{-\theta z} \phi\left(\frac{z}{\sqrt{r_n}}\right)  \frac{1}{(2k+1)!} \frac{1}{r_n^{(2k-1)/2}} H_{2k+1}\left(\frac{z}{\sqrt{r_n}}\right) \mathrm{d} z.
\end{align}
Now we show that we can drop the indicator function from the last line above.
Note that
for all $0\leq k\leq J$,
applying Lemma \ref{Properties-of-Hermite-Polynomial}(ii) and the inequality
\[
\int_A e^{-\theta z} \phi\left(\frac{z}{\sqrt{r_n}}\right) \left| H_{2k+1} \left(\frac{z}{\sqrt{r_n}}\right) \right| \mathrm{d} z\lesssim \int_A e^{-\theta z}\mathrm{d} z<\infty
\]
first and then Lemma \ref{General-many-to-one},
we get that
\begin{align}\label{step_17}
& 	\sum_{n=2}^\infty r_n^m \mathbb{E}_x\bigg( \sum_{v\in N(\sqrt{r_n})} \frac{1_{\{\min_{s\leq \sqrt{r_n}} X_v(s)>0\}}}{e^{(1-\frac{\theta^2}{2})\sqrt{r_n}}} e^{\theta X_v(\sqrt{r_n})}  \left| r_n^{(2k+1)/4}H_{2k+1}\left(\frac{X_v(\sqrt{r_n})}{r_n^{1/4}}\right) \right| \nonumber\\
	& \quad\quad \times 1_{\{ X_v(\sqrt{r_n})>  \sqrt{K\sqrt{r_n}\log n}\}}
	\int_A e^{-\theta z} \phi\left(\frac{z}{\sqrt{r_n}}\right)   \frac{1}{(2k+1)!} \frac{1}{r_n^{(2k-1)/2}} \left| H_{2k+1}\left(\frac{z}{\sqrt{r_n}}\right) \right| )\mathrm{d} z	\bigg)\nonumber\\
	&\lesssim \sum_{n=2}^\infty r_n^m \mathbb{E}_x\bigg( \sum_{v\in N(\sqrt{r_n})} \frac{1_{\{\min_{s\leq \sqrt{r_n}} X_v(s)>0\}}}{e^{(1-\frac{\theta^2}{2})\sqrt{r_n}}} e^{\theta X_v(\sqrt{r_n})}  r_n^{(2k+1)/4}\left(\left|\frac{X_v(\sqrt{r_n})}{r_n^{1/4}}\right|^{2k+1}+\left|\frac{X_v(\sqrt{r_n})}{r_n^{1/4}}\right|\right) \nonumber\\
	&\quad\quad\times   \frac{ 1}{r_n^{(2k-1)/2}} 1_{\{ X_v(\sqrt{r_n})>  \sqrt{K\sqrt{r_n}\log n}\}}
		\bigg)\nonumber\\
		& =  \sum_{n=2}^\infty \frac{r_n^m e^{\frac{\theta^2}{2}\sqrt{r_n}}}{r_n^{(2k-1)/2}} \mathbf{E}_x
		\bigg(1_{\left\{\min_{s\leq \sqrt{r_n}} X_\xi(s)>0 \right\}} X_\xi(\sqrt{r_n})e^{\theta X_\xi(\sqrt{r_n})} \nonumber\\
		&\quad\quad\times \left(X_\xi(\sqrt{r_n})^{2k}+ r_n^{k/2} \right)1_{\{X_\xi(\sqrt{r_n})>  \sqrt{K\sqrt{r_n}\log n}  \}}\bigg).
\end{align}
Recalling the choice of $\varepsilon$ in \eqref{Small-Varepsilon}.
 By Lemma \ref{lemma3},
 the right-hand side of \eqref{step_17} is equal to
\begin{align}\label{step_18}
	&\sum_{n=2}^\infty \frac{r_n^m }{r_n^{(2k-1)/2}} xe^{\theta x}\Pi_x^\uparrow \left(  \left( B_{\sqrt{r_n}}^{2k}+ r_n^{k/2}\right)1_{\{B_{\sqrt{r_n}}>  \sqrt{K\sqrt{r_n}\log n}  \}}\right)\nonumber\\
	& \leq \sum_{n=2}^\infty \frac{r_n^m }{r_n^{(2k-1)/2} n^{(1-\varepsilon)K/2} } xe^{\theta x} \Pi_x^\uparrow \left(  \left( B_{\sqrt{r_n}}^{2k}+ r_n^{k/2}\right)  e^{(1-\varepsilon) (B_{\sqrt{r_n}})^2/(2\sqrt{r_n})} \right)\nonumber\\
	& \lesssim  \sum_{n=2}^\infty \frac{r_n^m }{r_n^{(2k-1)/2} n^{(1-\varepsilon)K/2} } \int_0^\infty \frac{y^2}{r_n^{3/4}} e^{-(x-y)^2/(2\sqrt{r_n})}   \left( y^{2k}+ r_n^{k/2}\right)  e^{(1-\varepsilon) y^2/(2\sqrt{r_n})} \mathrm{d} y\nonumber\\
	& = \sum_{n=2}^\infty \frac{r_n^m r_n^{k/2} }{r_n^{(2k-1)/2} n^{(1-\varepsilon)K/2} } \int_0^\infty y^2 e^{-(xr_n^{-1/4}-y)^2/2}   \left( y^{2k}+ 1\right)  e^{(1-\varepsilon) y^2/2} \mathrm{d} y\nonumber\\
	& \leq \sum_{n=2}^\infty \frac{n^{(2m+1)/(2\kappa)}  }{ n^{(1-\varepsilon)K/2} } \int_0^\infty y^2 e^{-(xr_n^{-1/4}-y)^2/2}   \left( y^{2k}+ 1\right)  e^{(1-\varepsilon) y^2/2} \mathrm{d} y,
\end{align}
which is summable.
 Hence, combining \eqref{step_17} and \eqref{step_18},
we conclude that for all $0\leq k\leq J,$ almost surely,
\begin{align}
	& r_n^m \sum_{v\in N(\sqrt{r_n})} \frac{1_{\{\min_{s\leq \sqrt{r_n}} X_v(s)>0\}}}{e^{(1-\frac{\theta^2}{2})\sqrt{r_n}}} e^{\theta X_v(\sqrt{r_n})}  \left| r_n^{(2k+1)/4}H_{2k+1}\left(\frac{X_v(\sqrt{r_n})}{r_n^{1/4}}\right) \right| 1_{\{ X_v(\sqrt{r_n})>  \sqrt{K\sqrt{r_n}\log n}\}} \nonumber\\
	& \quad\quad \times
	\int_A e^{-\theta z}\phi\left(\frac{z}{\sqrt{r_n}}\right)   \frac{1}{(2k+1)!} \frac{1}{r_n^{(2k-1)/2}} \left| H_{2k+1}\left(\frac{z}{\sqrt{r_n}}\right) \right| \mathrm{d} z\stackrel{n\to\infty}{\longrightarrow}0.
\end{align}
Plugging this back to \eqref{step_19} and recalling the definition of the martingales in \eqref{Martingale}, we obtain that
\begin{align}\label{step_29}
		&\frac{Z_{r_n}^{(0,\infty)  } (A)}{r_n^{-3/2}e^{(1-\frac{\theta^2}{2})r_n}} =o(r_n^{-m})\nonumber\\
		&\quad +2 \sum_{k=0}^J M_{\sqrt{r_n}}^{(2k+1,\theta)}
		\int_A e^{-\theta z} \phi\left(\frac{z}{\sqrt{r_n}}\right)   \frac{1}{(2k+1)!} \frac{1}{r_n^{(2k-1)/2}} H_{2k+1}\left(\frac{z}{\sqrt{r_n}}\right)  \mathrm{d} z.
\end{align}

For all $m+1\leq k\leq J=2m+2\kappa +2$,
  we have $\lambda > J\geq k$. Therefore, for all $m+1\leq k\leq J$,
 by Lemma \ref{Properties-of-Hermite-Polynomial} (ii),
 Proposition
 \ref{Convergence-rate-martingale}(ii) and the fact that
    $\phi(y)|H_{2k+1}(y)| \lesssim \phi(y)|y| (|y|^{2k}+1)\lesssim |y|$,
\begin{align}\label{step_31}
	& \left| M_{\sqrt{r_n}}^{(2k+1,\theta)} \right|
	\int_A  e^{-\theta z} \phi\left(\frac{z}{\sqrt{r_n}}\right)   \frac{1}{(2k+1)!} \frac{1}{r_n^{(2k-1)/2}}  \left| H_{2k+1}\left(\frac{z}{\sqrt{r_n}}\right) \right| \mathrm{d} z\nonumber\\
	& \lesssim \left| M_{\sqrt{r_n}}^{(2k+1,\theta)} \right|	\int_A  e^{-\theta z}  \frac{1}{r_n^{(2k-1)/2}}  \frac{z}{\sqrt{r_n}}\mathrm{d} z = \frac{\left| M_{\sqrt{r_n}}^{(2k+1,\theta)} \right|}{r_n^k}	\int_A  z e^{-\theta z} \mathrm{d} z =o(r_n^{-m}).
\end{align}
Combining  \eqref{step_29} and \eqref{step_31}, we get that
\begin{align}\label{step_30}
	&\frac{Z_{r_n}^{(0,\infty)  } (A)}{r_n^{-3/2}e^{(1-\frac{\theta^2}{2})r_n}} =o(r_n^{-m})\nonumber\\
	&\quad + 2\sum_{k=0}^m M_{\sqrt{r_n}}^{(2k+1,\theta)}
	\int_A  e^{-\theta z} \phi\left(\frac{z}{\sqrt{r_n}}\right)   \frac{1}{(2k+1)!} \frac{1}{r_n^{(2k-1)/2}} H_{2k+1}\left(\frac{z}{\sqrt{r_n}}\right) \mathrm{d} z.
\end{align}
Noticing that $\lambda>2m$,
 let $\eta:= (\lambda -2m)/2<\lambda -2m$.
 By Proposition
\ref{Convergence-rate-martingale}(ii),
similar to \eqref{step_31},  for  all $0\leq k\leq m$,
\begin{align}\label{step_42}
	& \left| M_{\sqrt{r_n}}^{(2k+1,\theta)} -M_\infty^{(2k+1,\theta)} \right|
	\int_A  e^{-\theta z} \phi\left(\frac{z}{\sqrt{r_n}}\right)   \frac{1}{(2k+1)!} \frac{1}{r_n^{(2k-1)/2}}  \left| H_{2k+1}\left(\frac{z}{\sqrt{r_n}}\right) \right| \mathrm{d} z\nonumber\\
	& \lesssim  \frac{\left| M_{\sqrt{r_n}}^{(2k+1,\theta)} -M_\infty^{(2k+1,\theta)} \right|}{r_n^k}	\int_A z e^{-\theta z} \mathrm{d} z  = o(\sqrt{r_n}^{-(\lambda-k)+\eta} r_n^{-k})= o(r_n^{-m} r_n^{-\frac{k}{2}})= o(r_n^{-m}).
\end{align}
Since $\lambda> 2m+2\kappa +2$, we have
\begin{align}\label{step_32}
	&\frac{Z_{r_n}^{(0,\infty)  } (A)}{r_n^{-3/2}e^{(1-\frac{\theta^2}{2})r_n}} =o(r_n^{-m})\nonumber\\
	&\quad + 2\sum_{k=0}^m M_\infty^{(2k+1,\theta)}
		\int_A
	e^{-\theta z}\phi\left(\frac{z}{\sqrt{r_n}}\right)   \frac{1}{(2k+1)!} \frac{1}{r_n^{(2k-1)/2}} H_{2k+1}\left(\frac{z}{\sqrt{r_n}}\right)  \mathrm{d} z.
\end{align}
By \eqref{K-th-derivative} and Lemma \ref{Properties-of-Hermite-Polynomial}(i), for any $0\leq k\leq m$ and $x>0$, there exists $\xi\in (0,x)$ such that
\begin{align}
& \phi(x)H_{2k+1}(x)= \sum_{j=0}^{2m+1}\frac{\mathrm{d}^{j}}{\mathrm{d}x^j}\left(\phi H_{2k+1}\right)(0)\frac{x^j }{j!}+ \frac{x^{2m+2}}{(2m+2)!} \frac{\mathrm{d}^{2m+2}}{\mathrm{d}x^{2m+2}}\left(\phi H_{2k+1}\right)(\xi)\nonumber\\
& = (-1)^{2k+1}\sum_{j=0}^{m}\frac{\mathrm{d}^{2k+1+j}\phi(0)}{\mathrm{d}x^{2k+1+j}}\frac{x^{2j+1}}{(2j+1)!}+ O(x^{2m+2})= -\sum_{j=0}^m\frac{H_{2k+2j+2}(0)}{\sqrt{2\pi}}\frac{x^{2j+1}}{(2j+1)!}+ O(x^{2m+2}),
\end{align}
where in the second equality
we use the property that $H_{2\ell +1}(0)=0$.  Therefore,
\begin{align}
	&\frac{Z_{r_n}^{(0,\infty)  } (A)}{r_n^{-3/2}e^{(1-\frac{\theta^2}{2})r_n}} =o(r_n^{-m}) - \sum_{k=0}^m  \frac{2}{(2k+1)!} \frac{1}{r_n^{(2k-1)/2}}M_\infty^{(2k+1,\theta)} \nonumber\\
	&\quad\quad\times
	\int_A e^{-\theta z} \left(\sum_{j=0}^m\frac{H_{2k+2j+2}(0)}{(2j+1)!\sqrt{2\pi}}\frac{z^{2j+1}}{r_n^{(2j+1)/2}}+ O\left(\frac{z^{2m+2}}{r_n^{m+1}}\right)\right)\mathrm{d} z\nonumber\\
	& = o(r_n^{-m}) - \sqrt{\frac{2}{\pi}}\sum_{k=0}^m \sum_{j=0}^m \frac{1}{(2k+1)! (2j+1)!} \ \frac{H_{2k+2j+2}(0)}{r_n^{k+j}}M_\infty^{(2k+1,\theta)}
	\int_A  z^{2j+1} e^{-\theta z} \mathrm{d} z\nonumber\\
	& = o(r_n^{-m})  - \sqrt{\frac{2}{\pi}}\sum_{\ell=0}^m \frac{H_{2\ell+2}(0)}{r_n^{\ell}} \sum_{k=0}^{\ell} \frac{M_\infty^{(2k+1,\theta)} }{(2k+1)! (2\ell -2k+1)!}
	\int_A  z^{2\ell -2k+1} e^{-\theta z} \mathrm{d} z,
\end{align}
which completes the proof of (i) and (ii).

\textbf{Step 3:} In this step,
we prove (iii). Recall that when $\theta=0$ and $A=(a,\infty)$, $b_\theta=1/2$.
By \eqref{step_40} and Lemma \ref{lemma3} (i),
\begin{align}
	&\frac{Z_{r_n}^{(0,\infty)  } ((a,\infty))}{r_n^{-1/2}e^{r_n}} = o(r_n^{-m})+ \sum_{v\in N(\sqrt{r_n})} \frac{1_{\{\min_{s\leq \sqrt{r_n}} X_v(s) > 0 \}}}{r_n^{-1/2}e^{r_n}} 1_{\{ X_v(\sqrt{r_n})\leq  \sqrt{K\sqrt{r_n}\log n}\}}\nonumber\\
	&\quad \quad \times e^{r_n-\sqrt{r_n}}\Pi_{X_v(\sqrt{r_n})}\left( 1_{\{\min_{s\leq r_n -\sqrt{r_n}} B_s>0\}} 1_{(a,\infty)}(B_{r_n-\sqrt{r_n}}) \right)\nonumber\\
	& = o(r_n^{-m})+ \sum_{v\in N(\sqrt{r_n})} \frac{1_{\{\min_{s\leq \sqrt{r_n}} X_v(s) > 0 \}}}{r_n^{-1/2}e^{r_n}} 1_{\{ X_v(\sqrt{r_n})\leq  \sqrt{K\sqrt{r_n}\log n}\}}\nonumber\\
	&\quad \quad \times \frac{e^{r_n-\sqrt{r_n}}}{\sqrt{r_n -\sqrt{r_n}} }\int_a^\infty \left(\phi\left( \frac{y-X_v(\sqrt{r_n})}{\sqrt{r_n -\sqrt{r_n}}}\right) -\phi\left( \frac{y+X_v(\sqrt{r_n})}{\sqrt{r_n-\sqrt{r_n}}}\right)  \right)\mathrm{d}y\nonumber\\
	& = o(r_n^{-m})+ \sum_{v\in N(\sqrt{r_n})} \frac{1_{\{\min_{s\leq \sqrt{r_n}} X_v(s) > 0 \}}}{r_n^{-1/2}e^{\sqrt{r_n}}} 1_{\{ X_v(\sqrt{r_n})\leq  \sqrt{K\sqrt{r_n}\log n}\}}\nonumber\\
	&\quad \quad \times \left(\Phi\left( \frac{a+X_v(\sqrt{r_n})}{\sqrt{r_n -\sqrt{r_n}}}\right)-\Phi\left( \frac{a-X_v(\sqrt{r_n})}{\sqrt{r_n -\sqrt{r_n}}}\right) \right).
\end{align}
Put $J:= 2m+2\kappa+1> 2m+ \frac{K\kappa-1}{2}$.
By Lemma \ref{lemma4},  it holds that
\begin{align}
	&\frac{Z_{r_n}^{(0,\infty)  } ((a,\infty))}{r_n^{-1/2}e^{r_n}} = o(r_n^{-m})+ \sum_{v\in N(\sqrt{r_n})} \frac{1_{\{\min_{s\leq \sqrt{r_n}} X_v(s) > 0 \}}}{r_n^{-1/2}e^{\sqrt{r_n}}} 1_{\{ X_v(\sqrt{r_n})\leq  \sqrt{K\sqrt{r_n}\log n}\}} \nonumber\\
	&\quad  \times \left(
	2\phi\left(\frac{a}{\sqrt{r_n}}\right)  \sum_{k=0}^J  \frac{1}{(2k+1)!} \frac{1}{r_n^{(2k+1)/2}} H_{2k}\left(\frac{a}{\sqrt{r_n}}\right)  r_n^{(2k+1)/4}H_{2k+1}\left(\frac{X_u(\sqrt{r_n})}{r_n^{1/4}}   \right) +\varepsilon_{m,v, a, n} \right).
\end{align}
where the error term $\varepsilon_{m,v, a, n}$ satisfies that
\[
r_n^{(2m+1)/2} \sup_{a>0}\sup_{v\in N(\sqrt{r_n})} \left| \varepsilon_{m,v, a, n}\right|1_{\{ X_v(\sqrt{r_n})\leq  \sqrt{K\sqrt{r_n}\log n}\}}\stackrel{n\to\infty}{\longrightarrow}0,\quad
\P_x\mbox{-a.s.}
\]
Using the fact that $e^{-t}\sum_{v\in N(t)} 1$ is a non-negative martingale, we have
\begin{align}
	& r_n^m \sum_{v\in N(\sqrt{r_n})} \frac{1_{\{\min_{s\leq \sqrt{r_n}} X_v(s) > 0 \}}}{r_n^{-1/2}e^{\sqrt{r_n}}} 1_{\{ X_v(\sqrt{r_n})\leq  \sqrt{K\sqrt{r_n}\log n}\}}  \left|\varepsilon_{m,v, a, n}\right|\nonumber\\
	& \leq \left(e^{-\sqrt{r_n}}\sum_{v\in N(\sqrt{r_n})} 1 \right) r_n^{(2m+1)/2} \sup_{a>0}\sup_{v\in N(\sqrt{r_n})} \left| \varepsilon_{m,v, a, n}\right|1_{\{ X_v(\sqrt{r_n})\leq  \sqrt{K\sqrt{r_n}\log n}\}}\stackrel{n\to\infty}{\longrightarrow}0.
\end{align}
Therefore,
\begin{align}\label{step_41}
		&\frac{Z_{r_n}^{(0,\infty)  } ((a,\infty))}{r_n^{-1/2}e^{r_n}} = o(r_n^{-m})+ 2\sum_{v\in N(\sqrt{r_n})} \frac{1_{\{\min_{s\leq \sqrt{r_n}} X_v(s) > 0 \}}}{e^{\sqrt{r_n}}} 1_{\{ X_v(\sqrt{r_n})\leq  \sqrt{K\sqrt{r_n}\log n}\}} \nonumber\\
		&\quad  \times
		\phi\left(\frac{a}{\sqrt{r_n}}\right)  \sum_{k=0}^J  \frac{1}{(2k+1)!} \frac{1}{r_n^{k}} H_{2k}\left(\frac{a}{\sqrt{r_n}}\right)  r_n^{(2k+1)/4}H_{2k+1}\left(\frac{X_u(\sqrt{r_n})}{r_n^{1/4}}   \right)  .
\end{align}
Similar to the argument leading to \eqref{step_17}, for each $0\leq k\leq J$,
applying  Lemma \ref{Properties-of-Hermite-Polynomial} first, then Lemma \ref{General-many-to-one}  and Lemma \ref{lemma3}(i) at last, we also have that
\begin{align}
	&\sum_{n=1}^\infty r_n^m \E_x\bigg(\sum_{v\in N(\sqrt{r_n})} \frac{1_{\{\min_{s\leq \sqrt{r_n}} X_v(s) > 0 \}}}{e^{\sqrt{r_n}}} 1_{\{ X_v(\sqrt{r_n}) >  \sqrt{K\sqrt{r_n}\log n}\}} \nonumber\\
	&\quad  \times
	\phi\left(\frac{a}{\sqrt{r_n}}\right)  \frac{1}{r_n^{k}} \left| H_{2k}\left(\frac{a}{\sqrt{r_n}}\right) \right| \times \left| r_n^{(2k+1)/4}H_{2k+1}\left(\frac{X_u(\sqrt{r_n})}{r_n^{1/4}}   \right) \right|  \bigg)\nonumber\\
	& \lesssim \sum_{n=1}^\infty r_n^m \E_x\bigg(\sum_{v\in N(\sqrt{r_n})} \frac{1_{\{\min_{s\leq \sqrt{r_n}} X_v(s) > 0 \}}}{e^{\sqrt{r_n}}} 1_{\{ X_v(\sqrt{r_n})>  \sqrt{K\sqrt{r_n}\log n}\}} \nonumber\\
	&\quad  \times
	 \frac{1}{r_n^{k}}  X_u(\sqrt{r_n}) \left( \left(X_u(\sqrt{r_n})\right)^{2k}+ r_n^{k/2}  \right)   \bigg)\nonumber\\
	 & = \sum_{n=1}^\infty r_n^{m-k}
	 \mathbf{E}_x
	 \bigg( 1_{\{\min_{s\leq \sqrt{r_n}} X_\xi(s) > 0 \}} 1_{\{ X_\xi (\sqrt{r_n})>  \sqrt{K\sqrt{r_n}\log n}\}}  X_\xi(\sqrt{r_n}) \left( \left(X_\xi(\sqrt{r_n})\right)^{2k}+ r_n^{k/2}  \right)   \bigg)\nonumber\\
	 & = x\sum_{n=1}^\infty r_n^{m-k} \Pi_x^\uparrow \bigg( 1_{\{ B_{\sqrt{r_n}}>  \sqrt{K\sqrt{r_n}\log n}\}}   \left( \left(B_{\sqrt{r_n}}\right)^{2k}+ r_n^{k/2}  \right)   \bigg)<\infty,
\end{align}
where the last inequality follows from \eqref{step_18}.  Therefore, for $0\leq k\leq J$, $\P_x$-a.s.,
\begin{align}
	& r_n^m \sum_{v\in N(\sqrt{r_n})} \frac{1_{\{\min_{s\leq \sqrt{r_n}} X_v(s) > 0 \}}}{e^{\sqrt{r_n}}} 1_{\{ X_v(\sqrt{r_n}) >  \sqrt{K\sqrt{r_n}\log n}\}} \nonumber\\
	&\quad  \times
	\phi\left(\frac{a}{\sqrt{r_n}}\right)  \frac{1}{r_n^{k}} \left| H_{2k}\left(\frac{a}{\sqrt{r_n}}\right) \right| \times \left| r_n^{(2k+1)/4}H_{2k+1}\left(\frac{X_u(\sqrt{r_n})}{r_n^{1/4}}   \right) \right|   \stackrel{n\to\infty}{\longrightarrow}0.
\end{align}
Plugging this back to \eqref{step_41}, we finally get that
\begin{align}\label{step_43}
	&\frac{Z_{r_n}^{(0,\infty)  } ((a,\infty))}{r_n^{-1/2}e^{r_n}} = o(r_n^{-m})+2 \sum_{v\in N(\sqrt{r_n})} \frac{1_{\{\min_{s\leq \sqrt{r_n}} X_v(s) > 0 \}}}{e^{\sqrt{r_n}}}  \nonumber\\
	&\quad  \times
	\phi\left(\frac{a}{\sqrt{r_n}}\right)  \sum_{k=0}^J  \frac{1}{(2k+1)!} \frac{1}{r_n^{k}} H_{2k}\left(\frac{a}{\sqrt{r_n}}\right)  r_n^{(2k+1)/4}H_{2k+1}\left(\frac{X_u(\sqrt{r_n})}{r_n^{1/4}}   \right)  \nonumber\\
	& =o(r_n^{-m})+ 2\phi\left(\frac{a}{\sqrt{r_n}}\right)  \sum_{k=0}^J  \frac{1}{(2k+1)!} \frac{1}{r_n^{k}} H_{2k}\left(\frac{a}{\sqrt{r_n}}\right)  M_{\sqrt{r_n}}^{(2k+1,0)},
\end{align}
with $M_t^{(2k+1,0)}$ given in \eqref{Martingale}.  Using the  same argument as \eqref{step_31} and \eqref{step_42},
also noting that $\lambda>J= 2m+2\kappa+1$,
Proposition
\ref{Convergence-rate-martingale} (ii) and \eqref{step_43} imply that
\begin{align}
		&\frac{Z_{r_n}^{(0,\infty)  } ((a,\infty))}{r_n^{-1/2}e^{r_n}} = o(r_n^{-m})+  2\phi\left(\frac{a}{\sqrt{r_n}}\right)  \sum_{k=0}^m  \frac{1}{(2k+1)!} \frac{1}{r_n^{k}} H_{2k}\left(\frac{a}{\sqrt{r_n}}\right)  M_{\infty}^{(2k+1,0)}.
\end{align}
The remaining part is similar to the end of
Step 2
we omit the details here.
The proof of (iii) is complete.
\hfill$\Box$

\subsection{From discrete time to continuous time}

\begin{lemma}\label{Continuous-Lower-bound}
Let $x>0$ and $\theta \in [0, \sqrt{2})$.
For any given $m\in\N$, if $\kappa>2m+2$ and \eqref{LlogL} holds for
$\lambda>2m+2\kappa +2$,
then for any interval $J \subset (0,\infty)$, it holds that
	\begin{align}
		\liminf_{n\to\infty}  r_n^{m} \inf_{t\in (r_n, r_{n+1})}
		\frac{Z_{t}^{(0,\infty)} (J)- Z_{r_n}^{(0,\infty)} (J)}
		{r_n^{-3/2}e^{(1-\frac{\theta^2}{2})r_n}} \geq 0,\quad
		\P_x\mbox{-a.s.}
	\end{align}
\end{lemma}
\textbf{Proof: }
Suppose $x>0$, $m\in\N$, $\kappa>2m+2$ and that \eqref{LlogL} holds for
$\lambda>2m+2\kappa +2$.

\textbf{Step 1:} Define $\varepsilon_n:= \sqrt{r_{n+1}- r_n}$,
by the mean value theorem,
\begin{align}
\sqrt{r_n}\varepsilon_n =\sqrt{n^{1/\kappa}\left((n+1)^{1/\kappa}-n^{1/\kappa}\right)}\stackrel{\exists \xi\in [n,n+1]}{=}\sqrt{\frac{n^{1/\kappa}}{\kappa} \xi^{(-\kappa+1)/\kappa}  }\lesssim n^{(-\kappa+2)/(2\kappa)}\to 0.
\end{align}
For any $\eta< |J|/2$, define
\[
J_\eta:= \{y\in J: \mbox{dist}(y, J^c) \geq \eta \}.
\]
 For $u\in N(r_n)$, let $G_u$ be the event that
$u$ does not split before $r_{n+1}$ and that $\max_{s\in (r_n, r_{n+1})}|X_u(s)-X_u(r_n)|\le \sqrt{r_n}\varepsilon_n$. When $n$ is large enough so that
$\eta_n:= \sqrt{r_n}\varepsilon_n  < |J|/2$,
for $u\in N(r_n)$,
on the event $G_u$, for $t\in (r_n, r_{n+1})$, it must hold that
\[
\{X_u(r_n)\in J_{\eta_n} \}\subset    \{X_u(s) \in J,\  \forall s\in (r_n, r_{n+1})\}\subset \{X_u(t) \in J\} \cap \{ X_u(s)>0, \forall s\in (r_n, t]\}.
\]
Therefore, for $t\in (r_n, r_{n+1})$, by  the branching property,
\begin{align}\label{Lower-Continuous-time}
	& r_n^m \frac{Z_{t}^{(0,\infty)} (J)}{r_n^{-3/2}e^{(1-\frac{\theta^2}{2})r_n}} \nonumber\\
	&= \frac{r_n^{m+\frac{3}{2}}}{e^{(1-\frac{\theta^2}{2})r_n}} \sum_{u\in N(r_n)}1_{\{\min_{s\leq r_n} X_u(s)>0 \}} \sum_{v\in N(t): u\leq v} 1_{\{\min_{s< t-r_n} X_v(r_n+s)>0 \}}
	1_{\{X_v(t)\in J\}}\nonumber\\
	& \geq \frac{r_n^{m+\frac{3}{2}}}{e^{(1-\frac{\theta^2}{2})r_n}}  \sum_{u\in N(r_n)}1_{\{\min_{s\leq r_n} X_u(s)>0 \}}
	1_{\{X_u(r_n)\in J_{\eta_n}\}} 1_{G_u}=: I_n + II_n.
\end{align}
Here $I_n$ and $II_n$ are given by
\begin{align}
	I_n&:= \frac{r_n^{m+\frac{3}{2}}}{e^{(1-\frac{\theta^2}{2})r_n}}  \sum_{u\in N(r_n)}1_{\{\min_{s\leq r_n} X_u(s)>0 \}}
	1_{\{X_u(r_n)\in J_{\eta_n}\}}
	\left(1_{G_u}-\mathbb{P}_x\left(G_u\big| \mathcal{F}_{r_n} \right)\right),\nonumber\\
	II_n&:= \frac{r_n^{m+\frac{3}{2}}}{e^{(1-\frac{\theta^2}{2})r_n}}  \sum_{u\in N(r_n)}1_{\{\min_{s\leq r_n} X_u(s)>0 \}}	1_{\{X_u(r_n)\in J_{\eta_n}\}}	\mathbb{P}_x\left(G_u\big| \mathcal{F}_{r_n} \right).
\end{align}
We claim the following two limits hold:
\begin{align}\label{Claim-1}
	\sup_{t\in(r_n, r_{n+1})} |I_n| = |I_n|\stackrel{n\to\infty}{\longrightarrow}0,\quad
	\P_x\mbox{-a.s.}
\end{align}
and
\begin{align}\label{Claim-2}
	\sup_{t\in(r_n, r_{n+1})} \left|II_n -r_n^{m+\frac{3}{2}} \frac{Z_{r_n}^{(0,\infty)} (J)}{e^{(1-\frac{\theta^2}{2})r_n}}  \right| =\left|II_n -r_n^{m+\frac{3}{2}} \frac{Z_{r_n}^{(0,\infty)} (J)}{e^{(1-\frac{\theta^2}{2})r_n}}  \right|  \stackrel{n\to\infty}{\longrightarrow}0,\quad
	\P_x\mbox{-a.s.}
\end{align}
If \eqref{Claim-1} and \eqref{Claim-2} hold, then we complete the proof of Lemma together with \eqref{Lower-Continuous-time}.

\textbf{Step 2:} In this step, we prove \eqref{Claim-2}. Define
\[
III_n:= \frac{r_n^{m+\frac{3}{2}}}{e^{(1-\frac{\theta^2}{2})r_n}}  \sum_{u\in N(r_n)}1_{\{\min_{s\leq r_n} X_u(s)>0 \}}
1_{\{X_u(r_n)\in J\}}
 \mathbb{P}_x\left(G_u\big| \mathcal{F}_{r_n} \right).
\]
By Lemma \ref{lemma2}
with $A_n:= J\setminus J_{\eta_n}$,  as $n\to\infty$,
\begin{align}\label{step_36}
		\left|II_n-III_n \right| \leq  r_n^{m+\frac{3}{2}} 	\frac{Z_{r_n}^{(0,\infty)  } (A_n)}{e^{(1-\frac{\theta^2}{2})r_n}}= o(1)+ r_n^{m+\frac{3}{2}}	\frac{\mathbb{E}_x\left(Z_{r_n}^{(0,\infty)} (A_n)\Big| \mathcal{F}_{\sqrt{r_n}}\right)}{e^{(1-\frac{\theta^2}{2})r_n}},\quad
		\P_x\mbox{-a.s.}
\end{align}
By the Markov property and Lemma \ref{General-many-to-one}, we get that
\begin{align}\label{step_33}
	&r_n^{m+\frac{3}{2}}	\frac{\mathbb{E}_x\left(Z_{r_n}^{(0,\infty)} (A_n)\Big| \mathcal{F}_{\sqrt{r_n}}\right)}{e^{(1-\frac{\theta^2}{2})r_n}}
	= 	\frac{r_n^{m+\frac{3}{2}}  }{e^{(1-\frac{\theta^2}{2})r_n}} \sum_{u\in N(\sqrt{r_n})} 1_{\{\min_{s\leq \sqrt{r_n}} X_u(s)>0\}} \nonumber\\
	&\quad\quad\times e^{r_n-\sqrt{r_n}}	\mathbf{E}_{X_u(\sqrt{r_n})}
	\left(1_{\{\min_{s\leq r_n-\sqrt{r_n}} X_\xi(s)>0 \}} 1_{A_n}(X_\xi(r_n-\sqrt{r_n})) \right).
\end{align}
Since $(X_\xi(t), \mathbf{P}_x)$
is a standard Brownian motion with drift $-\theta$, by Lemma \ref{lemma3}(i) and (ii),
\begin{align}\label{step_34}
	& \mathbf{E}_{X_u(\sqrt{r_n})}
	\left(1_{\{\min_{s\leq r_n-\sqrt{r_n}} X_\xi(s)>0 \}} 1_{A_n}(X_\xi(r_n-\sqrt{r_n})) \right)\nonumber\\
	& = X_u(\sqrt{r_n})e^{\theta X_u(\sqrt{r_n})} e^{-\frac{\theta^2}{2}(r_n-\sqrt{r_n})}\Pi_{X_u(\sqrt{r_n})}^\uparrow \left( \frac{e^{-\theta B_{r_b-\sqrt{r_n}}}}{B_{r_n-\sqrt{r_n}}} 1_{A_n} (B_{r_n-\sqrt{r_n}})\right)\nonumber\\
	&= X_u(\sqrt{r_n})e^{\theta X_u(\sqrt{r_n})} e^{-\frac{\theta^2}{2}(r_n-\sqrt{r_n})} \int_{A_n} \frac{e^{-\theta y}}{y}
	 p_{r_n-\sqrt{r_n}}^\uparrow (X_u(\sqrt{r_n}), y)\mathrm{d}y\nonumber\\
	&\lesssim X_u(\sqrt{r_n})e^{\theta X_u(\sqrt{r_n})} e^{-\frac{\theta^2}{2}(r_n-\sqrt{r_n})} \int_{A_n} \frac{e^{-\theta y}}{y} \frac{y^2}{(r_n-\sqrt{r_n})^{3/2}}\mathrm{d}y.
\end{align}
Combining \eqref{step_33},  \eqref{step_34} and the definition of $M_t^{(1,\theta)}$ in \eqref{Martingale}, we see that
\begin{align}\label{step_35}
	& r_n^{m+\frac{3}{2}}	\frac{\mathbb{E}_x\left(Z_{r_n}^{(0,\infty)} (A_n)\Big| \mathcal{F}_{\sqrt{r_n}}\right)}{e^{(1-\frac{\theta^2}{2})r_n}}  \lesssim \frac{r_n^{m+\frac{3}{2}}  M_{\sqrt{r_n}}^{(1,\theta)}  }{ (r_n-\sqrt{r_n})^{3/2}} \int_{A_n} ye^{-\theta y} \mathrm{d}y\nonumber\\
	&\lesssim r_n^m |A_n| M_{\sqrt{r_n}}^{(1,\theta)}\lesssim r_n^m \sqrt{r_n}\varepsilon_n= \sqrt{n^{(2m+1)/\kappa} \left((n+1)^{1/\kappa}-n^{1/\kappa} \right) },
\end{align}
where in the second inequality,  we also used the fact that
$\sup_n \sup\{y: y\in A_n\}<\infty$
and that $\int_{A_n} y e^{-\theta y}\mathrm{d}y \leq  \sup_n \sup\{y: y\in A_n\}  \times |A_n|.$
Since $\kappa>2m+2$,
the last term of \eqref{step_35} tends to $0$. Combining \eqref{step_36}, \eqref{step_33} and \eqref{step_35}, it holds that
\begin{align}\label{step_37}
	|II_n- III_n| \stackrel{n\to\infty}{\longrightarrow}0,\quad
	\P_x\mbox{-a.s.}
\end{align}
By the Markov property and the definition of $G_u$,
\[
\mathbb{P}_x\left(G_u\big| \mathcal{F}_{r_n} \right)
= e^{-(r_{n+1}-r_n)} \Pi_0^{-\theta}\left(\max_{t<r_{n+1}-r_n}|B_t|\leq  \sqrt{r_n}\varepsilon_n\right).
\]
Thus, we get that
\begin{align}\label{step_38}
	&\left|III_n -r_n^{m+\frac{3}{2}}
	 \frac{Z_{r_n}^{(0,\infty)} (J)}{e^{(1-\frac{\theta^2}{2})r_n}}  \right| \nonumber\\
	&  \leq \frac{Z_{r_n}^{(0,\infty)} ((0,\infty))}{r_n^{-1/2}e^{(1-\frac{\theta^2}{2})r_n}} \cdot r_n^{m+1} \left(1- e^{-(r_{n+1}-r_n)} \Pi_0^{-\theta}\left(\max_{t<r_{n+1}-r_n}|B_t|\leq  \sqrt{r_n}\varepsilon_n\right)\right)\nonumber\\
	&\lesssim r_n^{m+1} \left(1- e^{-(r_{n+1}-r_n)} \Pi_0^{-\theta}\left(\max_{t<r_{n+1}-r_n}|B_t|\leq  \sqrt{r_n}\varepsilon_n\right)\right),
\end{align}
where in the last inequality we used Proposition \ref{Discrete-time-expansion} (i) (iii) with $m=0$.
Note that under the assumption $\kappa>2m+2$,
\begin{align}\label{step_39}
	 & r_n^{m+1}
	\left(1- e^{-(r_{n+1}-r_n)} \Pi_0^{-\theta}\left(\max_{t<r_{n+1}-r_n}|B_t|\leq  \sqrt{r_n}\varepsilon_n\right)\right)\nonumber\\
	&\leq  r_n^{m+1} O(r_{n+1}-r_n)+ r_n^{m+1}  \Pi_0^{-\theta}\left(\max_{t<r_{n+1}-r_n}|B_t|> \sqrt{r_n}\varepsilon_n\right)  \nonumber\\
	& \leq r_n^{m+1} O(n^{(-\kappa +1)/\kappa})+ r_n^{m+1} \Pi_0\left(\max_{t<r_{n+1}-r_n}|B_t|> \sqrt{r_n}\varepsilon_n - \theta (r_{n+1}-r_n)\right) \nonumber\\
	& = O(n^{-(\kappa-2-m)/\kappa})+ r_n^{m+1} \Pi_0\left(\max_{t<1}|B_t|> \sqrt{r_n} - \theta \sqrt{r_{n+1}-r_n}\right) = o (1).
\end{align}
Hence, \eqref{Claim-2} follows according to \eqref{step_37}, \eqref{step_38} and \eqref{step_39}.

\textbf{Step 3:} In this step, we prove \eqref{Claim-1}.
Since, given $\mathcal{F}_{r_n}$,
 $\left\{G_u : u\in N(r_n)\right\}$ are independent, we have
\begin{align}
	 &\mathbb{E}_x\left( \left|I_n\right|^2 \big|\mathcal{F}_{r_n}\right)=\frac{r_n^{2m+3}}{e^{2(1-\frac{\theta^2}{2})r_n}} \nonumber\\
	 &\quad\quad\times \sum_{u\in N(r_n)}1_{\{\min_{s\leq r_n} X_u(s)>0 \}}
	 1_{\{X_u(r_n)\in J_{\eta_n}\}}
	  \mathbb{E}_x\left( \left(1_{G_u}-\mathbb{P}_x\left(G_u\big| \mathcal{F}_n \right)\right)^2 \big|\mathcal{F}_{r_n}\right)\nonumber\\
	 & \leq 4 \frac{r_n^{2m+3}}{e^{2(1-\frac{\theta^2}{2})r_n}}  \sum_{u\in N(r_n)}1_{\{\min_{s\leq r_n} X_u(s)>0 \}} 1_{\{X_u(r_n)\in (0,\infty)\}}\nonumber\\
	 &= \frac{4 r_n^{2m+\frac{5}{2}}}{ e^{(1-\frac{\theta^2}{2})r_n}} \times  \frac{Z_{r_n}^{(0,\infty)} ((0,\infty))}{r_n^{-1/2}e^{(1-\frac{\theta^2}{2})r_n}}.
\end{align}
Now taking expectation with respect to $\P_x$, by Lemma \ref{General-many-to-one} and Lemma \ref{Sharp-Probability}(i)(ii), for any $\varepsilon>0$,
\begin{align}
	& \sum_{n=2}^\infty \mathbb{E}_x\left( \left|I_n\right|> \varepsilon \right)\leq \frac{1}{\varepsilon^2} \sum_{n=2}^\infty  \mathbb{E}_x\left( \left|I_n\right|^2 \right) \lesssim \sum_{n=2}^\infty   \frac{ r_n^{2m+\frac{5}{2}}}{ e^{(1-\frac{\theta^2}{2})r_n}} \frac{\mathbb{E}_x\left( Z_{r_n}^{(0,\infty)} ((0,\infty))\right)}{r_n^{-1/2}e^{(1-\frac{\theta^2}{2})r_n}}\nonumber\\
	&= \sum_{n=2}^\infty   \frac{ r_n^{2m+\frac{5}{2}}}{ e^{(1-\frac{\theta^2}{2})r_n}} \frac{e^{r_n}	\mathbf{P}_x \left(\min_{s\leq r_n} X_\xi(s)>0 \right)}{r_n^{-1/2}e^{(1-\frac{\theta^2}{2})r_n}}\lesssim \sum_{n=1}^\infty   \frac{ r_n^{2m+\frac{5}{2}}}{ e^{(1-\frac{\theta^2}{2})r_n}}<\infty,
\end{align}
which implies \eqref{Claim-1}. Thus, we complete the proof of the lemma.
\hfill$\Box$

\begin{lemma}\label{Continuous-time-approximation}
	Let $x>0$ and $\theta\in [0, \sqrt{2}).$ For any $m\in \N$,
	suppose that $\kappa>2m+2$ and \eqref{LlogL} holds with
	$\lambda>2m+2\kappa +2$,
	then for any interval $A\subset (0,\infty)$,
	$\P_x$-almost surely,
	\begin{align}
		r_n^{m} \sup_{t\in(r_n, r_{n+1})} \left| \frac{Z_{t}^{(0,\infty)} (A)- Z_{r_n}^{(0,\infty)} (A)}{r_n^{-3/2}e^{(1-\frac{\theta^2}{2})r_n}}\right|\stackrel{n\to\infty}{\longrightarrow} 0.
	\end{align}
\end{lemma}

\textbf{Proof: }
Suppose $x>0, m\in \N$, and that $\kappa>2m+2$ and \eqref{LlogL} holds with
$\lambda>2m+2\kappa +2$.
Fix $\theta \in [0,\sqrt{2})$.
Note that if
\[
\liminf_{n\to\infty} x_n \geq 0,\quad \liminf_{n\to\infty} y_n \geq 0,\quad \liminf_{n\to\infty} z_n \geq 0
\mbox{ and }\limsup_{n\to\infty} (x_n+y_n+z_n)\leq 0,
\]
then $\lim_{n\to\infty} x_n = \lim_{n\to\infty} y_n= \lim_{n\to\infty} z_n=0$.
Applying Lemma \ref{Continuous-Lower-bound} with $J=(a,b)$,
$[b,\infty)$ and $(0, a]$, we see that to prove  Lemma \ref{Continuous-time-approximation}, only need to prove that
\begin{align}\label{Continuous-Upper-bound}
		\limsup_{n\to\infty}  r_n^{m} \sup_{t\in (r_n, r_{n+1})}  \frac{Z_{t}^{(0,\infty)} ((0,\infty))- Z_{r_n}^{(0,\infty)} ((0,\infty))}{r_n^{-3/2}e^{(1-\frac{\theta^2}{2})r_n}} \leq 0.\quad \mbox{a.s.}
\end{align}
For any $t\in (r_{n}, r_{n+1})$, by the branching property, we see that
\begin{align}\label{Continuous-Upper-bound-2}
	& r_n^{m} \sup_{t\in (r_n, r_{n+1})}  \frac{Z_{t}^{(0,\infty)} ((0,\infty))- Z_{r_n}^{(0,\infty)} ((0,\infty))}{r_n^{-3/2}e^{(1-\frac{\theta^2}{2})r_n}} \nonumber\\
	&\leq \frac{r_n^{m+\frac{3}{2}}}{e^{(1-\frac{\theta^2}{2})r_n}}\sum_{u\in N(r_n)} 1_{\left\{\min_{s\leq r_n} X_u(s)>0\right\}} \sup_{t\in (r_n, r_{n+1})} \Big(\sum_{v\in N(t): u\leq v} 1 -1 \Big).
\end{align}
To drop the ``$\mbox{sup}$" above, we modify the branching particle system for
 $t\in (r_n, r_{n+1})$ such that when a particle dies in $(r_n,r_{n+1})$ and it splits into $L$ offspring, we modify the number of the offspring with $L+1$. For  $t\in(r_n, r_{n+1})$, we use  $\widetilde{N}(t)$ to denote the set of the particles alive at time $t$ in the modified process. It is obvious that the mean of the number in the  modified process is equal to $\sum_{k=0}^\infty (k+1)p_k=3$ and that for each $u\in N(r_n)$,
\[
\sup_{t\in (r_n, r_{n+1})} \Big(\sum_{v\in N(t): u\leq v} 1 -1 \Big) \leq \sup_{t\in (r_n, r_{n+1})} \Big(\sum_{v\in \widetilde{N}(t): u\leq v} 1 -1 \Big) =\Big( \sum_{v\in \widetilde{N}(r_{n+1}): u\leq v} 1 \Big)-1.
\]
Define
\begin{align}
		\widetilde{Z}_{n+1} &:= \frac{r_n^{m+\frac{3}{2}}}{e^{(1-\frac{\theta^2}{2})r_n}}\sum_{u\in N(r_n)} 1_{\left\{\min_{s\leq r_n} X_u(s)>0\right\}} \sum_{v\in \widetilde{N}(r_{n+1}): u\leq v} 1.
\end{align}
We claim that $\P_x$ almost surely,
\begin{align}\label{claim-tilde-Z}
	& \widetilde{Z}_{n+1} - \mathbb{E}_x\left(\widetilde{Z}_{n+1} \big|\mathcal{F}_{r_n}\right)  \to 0.
\end{align}
If the claim is true, then
\begin{align}
	& \widetilde{Z}_{n+1} - \mathbb{E}_x\left(\widetilde{Z}_{n+1} \big|\mathcal{F}_{r_n}\right)= \widetilde{Z}_{n+1}- e^{2(r_{n+1}-r_n)}r_n^{m}  \frac{Z_{r_n}^{(0,\infty)} ((0,\infty))}{r_n^{-3/2}e^{(1-\frac{\theta^2}{2})r_n}}   \to 0.
\end{align}
Using this and \eqref{Continuous-Upper-bound-2}, we get
\begin{align}
	& r_n^{m} \sup_{t\in (r_n, r_{n+1})}  \frac{Z_{t}^{(0,\infty)} ((0,\infty))- Z_{r_n}^{(0,\infty)} ((0,\infty))}{r_n^{-3/2}e^{(1-\frac{\theta^2}{2})r_n}} \leq \widetilde{Z}_{n+1}- r_n^{m+\frac{3}{2}}  \frac{Z_{r_n}^{(0,\infty)} ((0,\infty))}{e^{(1-\frac{\theta^2}{2})r_n}}  \nonumber\\ &
	=  \widetilde{Z}_{n+1} - \mathbb{E}_x\left(\widetilde{Z}_{n+1} \big|\mathcal{F}_{r_n}\right) + r_n^{m+1} (e^{2(r_{n+1}-r_n)}-1) \frac{Z_{r_n}^{(0,\infty)} ((0,\infty))}{r_n^{-1/2}e^{(1-\frac{\theta^2}{2})r_n}} \to 0,
\end{align}
where for the last limit, we used Proposition \ref{Discrete-time-expansion} (i) (iii) with $m=0$ and the fact that $r_n^{m+1}(e^{2(r_{n+1}-r_n)}-1) = r_n^{m+1} O(r_{n+1}-r_n)=o(1)$ under the assumption $\kappa>2m+2$.
Thus  the assertion of the lemma is valid.

Now we prove the claim \eqref{claim-tilde-Z}. We consider
another branching Brownian motion with underlying motion according to a standard Brownian motion with drift $-\theta$,  with  branching rate equal to $1$ and  with offspring distribution according to $\P_x(\widetilde{L}= k+1)=p_k$ for all $k\in \N$, then we may define another change-of-measure
\[
\frac{\mathrm{d}
\widetilde{\mathbf{P}}_x}{\mathrm{d}\P_x}
\bigg|_{\widetilde{\mathcal{F}}_t}:= \frac{\sum_{u\in \widetilde{N}(t)} 1}{e^{2t}},
\]
then a similar formula as Lemma \ref{General-many-to-one} can be established:

For any
$t>0$ and $u\in \widetilde{N}(t)$, let
$\Gamma(u,t)$
 be a  non-negative $\widetilde{\mathcal{F}}_t$-measurable  random variable. Then
\begin{align}\label{Another-Spine-decomposition}
		\E_x \Big(\sum_{u\in \widetilde{N}(t)} \Gamma(u,t) \Big) =e^{2t} \widetilde{\mathbf{E}}_x	\left(\Gamma (\xi_t, t)\right).
\end{align}
For $w\in \widetilde{N}(t)$ with $t\in (r_n, r_{n+1})$,  let $\widetilde{d}_w,\widetilde{O}_w$ denote the death time and the number of offspring of $w$ respectively.  For $v\in \widetilde{N}(r_{n+1})$, define
 $\widetilde{B}_{n,v}$ to be the event that, for all $w<v$ with $\widetilde{d}_w\in (r_n, r_{n+1})$, it holds that $\widetilde{O}_w \leq e^{c_0n}$, where $0< c_0< 1-\frac{\theta^2}{2}$ is fixed.
 Define \begin{align}
		\widetilde{Z}_{n+1}^{B} &:= \frac{r_n^{m+\frac{3}{2}}}{e^{(1-\frac{\theta^2}{2})r_n}}\sum_{u\in N(r_n)} 1_{\left\{\min_{s\leq r_n} X_u(s)>0\right\}} \sum_{v\in \widetilde{N}(r_{n+1}): u\leq v} 1_{\widetilde{B}_{n,v}}.
\end{align}
Now for $v\in \widetilde{N}(r_{n+1}-r_n)$, $\widetilde{D}_{n,v}$ denotes the event that, for all $w<v$, it holds that $\widetilde{O}_w\leq e^{c_0 n}$.
Let $\widetilde{d}_i$ be the $i$-th splitting time of the spine and $\widetilde{O}_i$ be the number of children produced by the spine at time $\widetilde{d}_i$.
Define $\widetilde{D}_{n, \xi_{r_{n+1}-r_n}}$ to be the event that, for all $i$ with $\widetilde{d}_i<r_{n+1}-r_n$, it holds that $\widetilde{O}_i\leq e^{c_0 n}$.
Then by the  branching property, the Markov property and \eqref{Another-Spine-decomposition}, we have
\begin{align}
	&\mathbb{E}_x \left( \left| \widetilde{Z}_{n+1} -\widetilde{Z}_{n+1}^B  \right|  \big| \mathcal{F}_{r_n}\right)
	\nonumber\\
	& = \frac{r_n^{m+\frac{3}{2}}}{e^{(1-\frac{\theta^2}{2})r_n}}\sum_{u\in N(r_n)} 1_{\left\{\min_{s\leq r_n} X_u(s)>0\right\}} \E_{X_u(r_n)}\Big(\sum_{v\in \widetilde{N}(r_{n+1}-r_n)} 1_{\widetilde{D}_{n,v}^c}\Big)\nonumber\\
	& = e^{2(r_{n+1}-r_n)}\frac{r_n^{m+\frac{3}{2}}}{e^{(1-\frac{\theta^2}{2})r_n}}\sum_{u\in N(r_n)} 1_{\left\{\min_{s\leq r_n} X_u(s)>0\right\}}
	\widetilde{\mathbf{E}}_{X_u(r_n)}
	\Big(1_{\widetilde{D}_{n,\xi_{r_{n+1}-r_n}}^c}\Big).
\end{align}
Noticing that \eqref{LlogL} implies $\widetilde{\mathbf{E}}_x(\log_+^{1+\lambda} \widetilde{O}_1) = \widetilde{\mathbf{E}}_0(\log_+^{1+\lambda} \widetilde{O}_1) <\infty$, we obtain
\begin{align}
	&\mathbb{E}_x \left( \left| \widetilde{Z}_{n+1} -\widetilde{Z}_{n+1}^B  \right|  \big| \mathcal{F}_{r_n}\right)
	\nonumber\\
	&\lesssim \frac{r_n^{m+\frac{3}{2}}}{e^{(1-\frac{\theta^2}{2})r_n}}\sum_{u\in N(r_n)} 1_{\left\{\min_{s\leq r_n} X_u(s)>0\right\}}	\widetilde{\mathbf{E}}_{X_u(r_n)	}\Big(\sum_{i: \widetilde{d}_i <r_{n+1}-r_n} 1_{\{\widetilde{O}_i > e^{c_0n}\}} \Big)\nonumber\\
	&\lesssim \frac{r_n^{m+\frac{3}{2}}}{e^{(1-\frac{\theta^2}{2})r_n}}\sum_{u\in N(r_n)} 1_{\left\{\min_{s\leq r_n} X_u(s)>0\right\}}  \frac{r_{n+1}-r_n}{n^{1+\lambda}}.
\end{align}
Now taking expectation with respect to $\P_x$, using Lemma \ref{General-many-to-one} and Lemma \ref{Sharp-Probability}(i) (ii), we get that
\begin{align}\label{First-moment-2}
	&\mathbb{E}_x \left( \left| \widetilde{Z}_{n+1} -\widetilde{Z}_{n+1}^B  \right|  \right)
	\lesssim  \frac{r_{n+1}-r_n}{n^{1+\lambda}}\frac{r_n^{m+\frac{3}{2}}}{e^{(1-\frac{\theta^2}{2})r_n}} e^{r_n}
	\mathbf{P}_x
	\left(\min_{s\leq r_n} X_\xi (s)>0\right)\nonumber\\
	&\lesssim xe^{\theta x} \frac{r_{n+1}-r_n}{n^{1+\lambda}} r_n^{m+1}.
\end{align}
Similarly, by the branching property and the Markov property,
\begin{align}\label{Second-moment-1}
	&\mathbb{E}_x \left( \left| \widetilde{Z}_{n+1}^B - \mathbb{E}_x\left(\widetilde{Z}_{n+1}^B \big|\mathcal{F}_{r_n}\right) \right|^2  \big| \mathcal{F}_{r_n}\right)= \frac{r_n^{2m+3}}{e^{2(1-\frac{\theta^2}{2})r_n}}\sum_{u\in N(r_n)} 1_{\left\{\min_{s\leq r_n} X_u(s)>0\right\}} \nonumber\\
	&\quad\quad\times \mathbb{E}_{X_u(r_n)}\Big( \Big(\sum_{v\in \widetilde{N}(r_{n+1}-r_n)} 1_{\widetilde{D}_{n,v}}- \E_{X_u(r_n)}\Big(\sum_{v\in \widetilde{N}(r_{n+1}-r_n)} 1_{\widetilde{D}_{n,v}}\Big) \Big)^2\Big)\nonumber\\
	&\leq \frac{r_n^{2m+3}}{e^{2(1-\frac{\theta^2}{2})r_n}}\sum_{u\in N(r_n)} 1_{\left\{\min_{s\leq r_n} X_u(s)>0\right\}} \mathbb{E}_{X_u(r_n)}\Big( \Big(\sum_{v\in \widetilde{N}(r_{n+1}-r_n)} 1_{\widetilde{D}_{n,v}} \Big)^2\Big).
\end{align}
By \eqref{Another-Spine-decomposition}, we see that
\begin{align}
	& \mathbb{E}_{X_u(r_n)}\Big( \Big(\sum_{v\in \widetilde{N}(r_{n+1}-r_n)} 1_{\widetilde{D}_{n,v}} \Big)^2\Big)= e^{2(r_{n+1}-r_n)}
	\widetilde{\mathbf{E}}_{X_u(r_n)}
	\Big( 1_{\widetilde{D}_{n,\xi_{r_{n+1}-r_n} }}  \sum_{v\in \widetilde{N}(r_{n+1}-r_n)} 1_{\widetilde{D}_{n,v}}\Big)\nonumber\\
	&\lesssim
	\widetilde{\mathbf{E}}_{X_u(r_n)}
	\Big( 1_{\widetilde{D}_{n,\xi_{r_{n+1}-r_n} }}  \sum_{v\in \widetilde{N}(r_{n+1}-r_n)} 1 \Big).
\end{align}
On the set $\widetilde{D}_{n,\xi_{r_{n+1}-r_n} }$, we have
\begin{align}
	&\widetilde{\mathbf{E}}_{X_u(r_n)}
	\Big(\sum_{v\in \widetilde{N}(r_{n+1}-r_n)} 1  \big| \widetilde{d}_i, \widetilde{O}_i: i\geq 1\Big)=\sum_{i: \widetilde{d}_i< r_{n+1}-r_n} \left(\widetilde{O}_i -1\right) e^{2(r_{n+1}-r_n-d_i)}\nonumber\\
	& \lesssim e^{c_0 n} \sum_{i: \widetilde{d}_i< r_{n+1}-r_n} 1,
\end{align}
which implies that
\begin{align}
	 \widetilde{\mathbf{E}}_{X_u(r_n)}
	 \Big( 1_{\widetilde{D}_{n,\xi_{r_{n+1}-r_n} }}  \sum_{v\in \widetilde{N}(r_{n+1}-r_n)} 1 \Big)\lesssim e^{c_0 n}  	 \widetilde{\mathbf{E}}_{X_u(r_n)}
	 \Big(\sum_{i: \widetilde{d}_i< r_{n+1}-r_n} 1\Big)\lesssim e^{c_0 n}(r_{n+1}-r_n).
\end{align}
Therefore, plugging this upper-bound back to \eqref{Second-moment-1}, we have
\begin{align}
	&\mathbb{E}_x \left( \left| \widetilde{Z}_{n+1}^B - \mathbb{E}_x\left(\widetilde{Z}_{n+1}^B \big|\mathcal{F}_{r_n}\right) \right|^2  \big| \mathcal{F}_{r_n}\right)\lesssim \frac{r_n^{2m+3}e^{c_0 n}(r_{n+1}-r_n)}{e^{2(1-\frac{\theta^2}{2})r_n}}\sum_{u\in N(r_n)} 1_{\left\{\min_{s\leq r_n} X_u(s)>0\right\}}.
\end{align}
Taking expectation with respect to $\P_x$, using Lemma \ref{General-many-to-one} and Lemma \ref{Sharp-Probability}(i), we conclude that
\begin{align}\label{Second-moment-2}
	&\mathbb{E}_x \left( \left| \widetilde{Z}_{n+1}^B - \mathbb{E}_x\left(\widetilde{Z}_{n+1}^B \big|\mathcal{F}_{r_n}\right) \right|^2 \right) \lesssim  \frac{r_n^{2m+3}e^{c_0 n}(r_{n+1}-r_n)}{e^{2(1-\frac{\theta^2}{2})r_n}}   e^{r_n}	\mathbf{E}_x	 \left( 1_{\left\{\min_{s\leq r_n} X_\xi(s)>0\right\}}\right)\nonumber\\
	&\lesssim  \frac{r_n^{2m+3}e^{c_0 n}(r_{n+1}-r_n)}{e^{(1-\frac{\theta^2}{2})r_n}}  .
\end{align}
First using \eqref{New-step_2} with $X= \widetilde{Z}_{n+1} , Y= \widetilde{Z}_{n+1}^B$ and $\ell=2$, and then \eqref{First-moment-2} and \eqref{Second-moment-2}, we get that
\begin{align}
	& \sum_{n=1}^\infty \mathbb{E}_x \left( \left| \widetilde{Z}_{n+1} - \mathbb{E}_x\left(\widetilde{Z}_{n+1} \big|\mathcal{F}_{r_n}\right) \right| \right)  \nonumber\\
	& \leq 2\sum_{n=1}^\infty \mathbb{E}_x \left( \left| \widetilde{Z}_{n+1} -\widetilde{Z}_{n+1}^B  \right|  \right) + \sum_{n=1}^\infty \mathbb{E}_x \left( \left| \widetilde{Z}_{n+1}^B - \mathbb{E}_x\left(\widetilde{Z}_{n+1}^B \big|\mathcal{F}_{r_n}\right) \right|^2  \right)^{1/2}\nonumber\\
	& \lesssim  \sum_{n=1}^\infty  \left(\frac{r_{n+1}-r_n}{n^{1+\lambda}} r_n^{m+1} + \left( \frac{r_n^{2m+3}e^{c_0 n}(r_{n+1}-r_n)}{e^{(1-\frac{\theta^2}{2})r_n}}   \right)^{1/2}\right),
\end{align}
which is finite since for $\lambda>0$ and $\kappa >2m+2$, $(r_{n+1}-r_n)r_n^{m+1} = O(n^{(m-\kappa+2)/\kappa})=o(1)$ and $\sum_{n=1}^\infty \frac{1}{n^{1+\lambda}}<\infty$.
Then we finished the proof of the claim \eqref{claim-tilde-Z}.
\hfill$\Box$

\textbf{Proof of Theorem \ref{thm1}:}
For any $\lambda>6m+6$, we can find
an appropriate
$\kappa$ satisfying the conditions
of Proposition \ref{Discrete-time-expansion} and Lemma \ref{Continuous-time-approximation}. For instance,
we can take $\kappa:= 2m+2+ (\lambda-6m-6)/4>2m+2$, then
\begin{align}
		\lambda = 4(\kappa-2m-2) +6m+6 > 6m+6+2(\kappa-2m-2)=2m+2\kappa +2.
\end{align}
For any $\ell\in [0,m]$, since $\kappa>2m+2$, we have
\begin{equation}\label{difference}
 r_{n+1}^\ell-r_n^\ell=\frac{\ell}\kappa\int^{n+1}_ny^{-(\kappa-\ell)/\kappa}dy\le \frac{m}\kappa n^{-(\kappa-m)/\kappa}
 = o(r_n^{-m}).
 \end{equation}
For $t\in (r_{n}, r_{n+1})$,  by Lemma \ref{Continuous-time-approximation} and \eqref{difference},
we get that for any $a\ge 0$ and $t\in (r_{n}, r_{n+1})$,
\begin{align}
	&\frac{Z_{t}^{(0,\infty)} ((a, \infty))}{t^{-3/2}e^{(1-\frac{\theta^2}{2})t}}= \frac{r_n^{-3/2}}{t^{-3/2}}\cdot \frac{e^{(1-\frac{\theta^2}{2})r_n}}{e^{(1-\frac{\theta^2}{2})t}} \cdot \frac{Z_{t}^{(0,\infty)} ((a, \infty))}{r_n^{-3/2}e^{(1-\frac{\theta^2}{2})r_n}}\nonumber\\
	&=	\left(1+o (r_n^{-m})\right)	\cdot \left(o\left(r_n^{-m}\right) + \frac{Z_{r_n}^{(0,\infty)} ((a, \infty))}{r_n^{-3/2}e^{(1-\frac{\theta^2}{2})r_n}} \right)\nonumber\\
	& = o(r_n^{-m})+ \frac{Z_{r_n}^{(0,\infty)} ((a, \infty))}{r_n^{-3/2}e^{(1-\frac{\theta^2}{2})r_n}}.
\end{align}
By \eqref{difference}, we have
\[
r_n^m \sup_{t\in(r_n, r_{n+1})} \left| \frac{1}{t^\ell} -\frac{1}{r_n^\ell} \right|
= o(1).
\]
Combining the above with Proposition  \ref{Discrete-time-expansion}(i), we get the assertion of Theorem \ref{thm1}.
\hfill$\Box$

\textbf{Proof of Theorem \ref{thrm2}:}  The proof is similar to that of Theorem \ref{thm1} and we omit the
details.
\hfill$\Box$

\bigskip
\noindent
{\bf Acknowledgements:}
We thank Zhenyao Sun for helpful discussions.
\bigskip
\noindent

\begin{singlespace}

\end{singlespace}

\vskip 0.2truein
\vskip 0.2truein

\noindent{\bf Haojie Hou:}  School of Mathematical Sciences, Peking
University,   Beijing, 100871, P.R. China. Email: {\texttt
houhaojie@pku.edu.cn}

\smallskip

\noindent{\bf Yan-Xia Ren:} LMAM School of Mathematical Sciences \& Center for
Statistical Science, Peking
University,  Beijing, 100871, P.R. China. Email: {\texttt
yxren@math.pku.edu.cn}

\smallskip
\noindent {\bf Renming Song:} Department of Mathematics,
University of Illinois Urbana-Champaign,
Urbana, IL 61801, U.S.A.
Email: {\texttt rsong@illinois.edu}

\end{document}